%
\catcode`@=11
%
%
\def\bibn@me{R\'ef\'erences}
\def\bibliographym@rk{\centerline{{\sc\bibn@me}}
	\sectionmark\section{\ignorespaces}{\unskip\bibn@me}
	\bigbreak\bgroup
	\ifx\ninepoint\undefined\relax\else\ninepoint\fi}
%
%
%
\let\refsp@ce=\ 
\let\bibleftm@rk=[
\let\bibrightm@rk=]
%
%
%
\def\numero{n\raise.82ex\hbox{$\fam0\scriptscriptstyle o$}~\ignorespaces}
%
%
\newcount\equationc@unt
\newcount\bibc@unt
\newif\ifref@changes\ref@changesfalse
\newif\ifpageref@changes\ref@changesfalse
\newif\ifbib@changes\bib@changesfalse
\newif\ifref@undefined\ref@undefinedfalse
\newif\ifpageref@undefined\ref@undefinedfalse
\newif\ifbib@undefined\bib@undefinedfalse
\newwrite\@auxout
%
%
\def\eqnum{\global\advance\equationc@unt by 1%
\edef\lastref{\number\equationc@unt}%
\eqno{(\lastref)}}
%
%
%
%
%
%
\def\re@dreferences#1#2{{%
	\re@dreferenceslist{#1}#2,\undefined\@@}}
\def\re@dreferenceslist#1#2,#3\@@{\def\next{#2}%
	\expandafter\ifx\csname#1@@\meaning\next\endcsname\relax
	??\immediate\write16
	{Warning, #1-reference "\next" on page \the\pageno\space
	is undefined.}%
	\global\csname#1@undefinedtrue\endcsname
	\else\csname#1@@\meaning\next\endcsname\fi
	\ifx#3\undefined\relax
	\else,\refsp@ce\re@dreferenceslist{#1}#3\@@\fi}
%
%
%
\def\newlabel#1#2{{\def\next{#1}\newl@bel#2}}
\def\newl@bel#1#2{%
	\expandafter\xdef\csname ref@@\meaning\next\endcsname{#1}%
	\expandafter\xdef\csname pageref@@\meaning\next\endcsname{#2}}
\def\label#1{{%
	\toks0={#1}\message{ref(\lastref) \the\toks0,}%
	\ignorespaces\immediate\write\@auxout%
	{\noexpand\newlabel{\the\toks0}{{\lastref}{\the\pageno}}}%
	\def\next{#1}%
	\expandafter\ifx\csname ref@@\meaning\next\endcsname\lastref%
	\else\global\ref@changestrue\fi%
	\newlabel{#1}{{\lastref}{\the\pageno}}}}
\def\ref#1{\re@dreferences{ref}{#1}}
\def\pageref#1{\re@dreferences{pageref}{#1}}
%
%
\def\bibcite#1#2{{\def\next{#1}%
	\expandafter\xdef\csname bib@@\meaning\next\endcsname{#2}}}
\def\cite#1{\bibleftm@rk\re@dreferences{bib}{#1}\bibrightm@rk}
%
%
\def\beginthebibliography#1{\bibliographym@rk
	\setbox0\hbox{\bibleftm@rk#1\bibrightm@rk\enspace}
	\parindent=\wd0
	\global\bibc@unt=0
	\def\bibitem##1{\global\advance\bibc@unt by 1
		\edef\lastref{\number\bibc@unt}
		{\toks0={##1}
		\message{bib[\lastref] \the\toks0,}%
		\immediate\write\@auxout
		{\noexpand\bibcite{\the\toks0}{\lastref}}}
		\def\next{##1}%
		\expandafter\ifx
		\csname bib@@\meaning\next\endcsname\lastref
		\else\global\bib@changestrue\fi%
		\bibcite{##1}{\lastref}
		\medbreak
		\item{\hfill\bibleftm@rk\lastref\bibrightm@rk}%
		}
	}
\def\endthebibliography{\egroup\par}
%
%
%
\def\@closeaux{\closeout\@auxout
	\ifref@changes\immediate\write16
	{Warning, changes in references.}\fi
	\ifpageref@changes\immediate\write16
	{Warning, changes in page references.}\fi
	\ifbib@changes\immediate\write16
	{Warning, changes in bibliography.}\fi
	\ifref@undefined\immediate\write16
	{Warning, references undefined.}\fi
	\ifpageref@undefined\immediate\write16
	{Warning, page references undefined.}\fi
	\ifbib@undefined\immediate\write16
	{Warning, citations undefined.}\fi}
%
%
\immediate\openin\@auxout=\jobname.aux
\ifeof\@auxout \immediate\write16
  {Creating file \jobname.aux}
\immediate\closein\@auxout
\immediate\openout\@auxout=\jobname.aux
\immediate\write\@auxout {\relax}%
\immediate\closeout\@auxout
\else\immediate\closein\@auxout\fi
%
%
\input\jobname.aux
\immediate\openout\@auxout=\jobname.aux
%
%
\catcode`@=12

\def\Q{{\bf {Q}}}

\def\Z{{\bf Z}} 
\def\R{{\bf R}} 
\def\C{{\bf C}}  

\def\Re{{\rm Re}\,}
\def\Im{{\rm Im}\,}
\def\variable{X}

\def\vep{{\varepsilon}}

\catcode`@=11
\def\bibliographym@rk{\bgroup}


\outer\def\bye{     \par\vfill\supereject\end}

\def\og{\leavevmode\raise.3ex\hbox{$\scriptscriptstyle
\langle\!\langle\,$}}
\def \fg {\leavevmode\raise.3ex\hbox{$\scriptscriptstyle
\!\rangle\!\rangle\,\,$}}

\magnification=1200
\parskip 4pt

\frenchspacing

\def\house#1{\setbox1=\hbox{$\,#1\,$}%
\dimen1=\ht1 \advance\dimen1 by 2pt \dimen2=\dp1 \advance\dimen2 by 2pt
\setbox1=\hbox{\vrule height\dimen1 depth\dimen2\box1\vrule}%
\setbox1=\vbox{\hrule\box1}%
\advance\dimen1 by .4pt \ht1=\dimen1
\advance\dimen2 by .4pt \dp1=\dimen2 \box1\relax}

\def\hw{{\hat{w}}}
\def\tw{{\widetilde{w}}}
\def\what{{\hat{w}}}

\def\eps{{\varepsilon}}
\def\vep{{\varepsilon}}

\def\sm{\smallskip}  \def\noi{\noindent}

\def\build#1_#2^#3{\mathrel{\mathop{\kern 0pt#1}\limits_{#2}^{#3}}}

\font\fivegoth=eufm5 \font\sevengoth=eufm7 \font\tengoth=eufm10

\newfam\gothfam \scriptscriptfont\gothfam=\fivegoth
\textfont\gothfam=\tengoth \scriptfont\gothfam=\sevengoth

\def\ovxi{{\overline \xi}}

\def\proof{\noindent {\bf Proof: }}

\def\smallsquare{\vbox{\hrule\hbox{\vrule height 1 ex\kern 1 ex\vrule}\hrule}}
\def\cqfd{\hfill \smallsquare\vskip 3mm}
\def\qed{\hfill \smallsquare\vskip 3mm}

\vskip 2mm

\centerline{\bf Approximation of complex algebraic numbers}

\sm

\centerline
{\bf by algebraic numbers of bounded degree}

\vskip 8mm

\centerline{Y{\sevenrm ANN} B{\sevenrm UGEAUD} (Strasbourg) \ {\it \&}
\ J{\sevenrm AN}-H{\sevenrm ENDRIK} E{\sevenrm VERTSE} (Leiden) \footnote{}{\rm
}}

{\narrower\narrower
\vskip 11mm

\proclaim Abstract.
{\rm
To measure how well a given complex number $\xi$ can be approximated
by algebraic numbers of degree at most $n$ one may use
the quantities $w_n(\xi )$ and $w^*_n(\xi )$ introduced by Mahler
and Koksma, respectively.
The values of $w_n(\xi )$ and $w_n^*(\xi)$ have been computed for
real algebraic numbers $\xi$, 
but up to now not for complex,
non-real algebraic
numbers $\xi$. In this paper we compute $w_n(\xi )$, $w_n^*(\xi )$
for all positive integers $n$ and algebraic numbers $\xi\in\C\setminus\R$,
except for those pairs $(n,\xi )$ such that
$n$ is even, $n\geq 6$ and $n+3\le \deg\xi \le 2n-2$.
It is known that every real algebraic number of degree $>n$ has the same
values for $w_n$ and $w_n^*$ as almost every real number.  
Our results imply that for every
positive even integer $n$ there are complex algebraic numbers
$\xi$ of degree $>n$
which are unusually well approximable by algebraic numbers of degree
at most $n$, i.e., have larger values for $w_n$ and $w_n^*$ than almost
all complex numbers.
We consider also the approximation of complex non-real algebraic numbers $\xi$
by algebraic integers, and show that if $\xi$ is unusually well approximable
by algebraic numbers of degree at most $n$ then
it is unusually badly approximable by algebraic integers
of degree at most $n+1$.
By means of Schmidt's Subspace Theorem we reduce the approximation
problem to compute $w_n(\xi )$, $w_n^*(\xi )$ to an algebraic problem
which is trivial if $\xi$ is real but much harder if $\xi$ is not real.
We give a partial solution to this problem.
}

}

\vskip 6mm

\centerline{\bf 1. Introduction}

\vskip 5mm

Conjecturally, most of the properties shared by almost
all numbers (throughout the present paper, `almost all'
always refers to the Lebesgue measure) should be 
either trivially false for the algebraic numbers, or satisfied by the
algebraic numbers. Thus, the sequence of partial quotients of
every real, irrational algebraic number
of degree at least $3$ 
is expected to be unbounded,
and the digit $2$ should occur infinitely often in the decimal expansion
of every real, irrational algebraic number. Our very limited knowledge on
these two problems show that they are far from being solved.

In Diophantine approximation, the situation is better understood.
For instance, for $\xi\in\R$, denote by $\lambda (\xi )$ the supremum
of all $\lambda$ such that the inequality
$|\xi -p/q |\leq \max \{|p|,|q|\}^{-\lambda}$ has
infinitely many solutions in rational numbers $p/q$ where $p,q\in\Z$, $q\not=0$.
Then for almost all real numbers $\xi$
we have $\lambda (\xi )=2$,
while by Roth's Theorem \cite{Ro}, we have also $\lambda (\xi )=2$
for every real, algebraic, irrational number $\xi$.

More generally, the quality of the approximation
of a complex number $\xi$ by algebraic numbers of degree at most
$n$ can be measured by means of the exponents $w_n(\xi )$ and $w_n^*(\xi )$
introduced by Mahler \cite{Mah} in 1932 and
by Koksma \cite{Kok} in 1939, respectively, which are defined as follows:
\vskip 0.2cm\noindent
$\bullet$
$w_n(\xi)$ denotes the supremum of those real numbers $w$ for
which the inequality
$$
0 < |P(\xi)| \le H(P)^{-w}
$$
is satisfied by infinitely many polynomials
$P\in\Z [X]$ of degree at most $n$;
\vskip 0.1cm\noindent
$\bullet$
$w_n^*(\xi)$ denotes 
the supremum of those real numbers $w^*$
for which the inequality
$$
0 < |\xi - \alpha| \le H(\alpha)^{- w^* -1}
$$
is satisfied by infinitely many algebraic numbers $\alpha$
of degree at most $n$.
\vskip0.2cm\noindent
Here, the height $H(P)$ of a polynomial $P\in\Z [X]$
is defined to be the maximum of the absolute values of its coefficients,
and the height $H(\alpha )$ of an algebraic number $\alpha$
is defined to be the height of its minimal polynomial (by definition
with coprime integer coefficients).
The reader is directed to \cite{BuLiv} for an overview of the known results
on the functions $w_n$ and $w_n^*$.

For every complex number $\xi$ and every integer $n \geq 1$
one has $w_n^*(\xi) \leq w_n(\xi)$,
but for every $n \geq 2$, there are complex numbers $\xi$
for which the inequality is strict.
Sprind\v zuk (see his monograph \cite{Spr})
established in 1965 that for every integer $n \geq 1$, we have
$w_n (\xi) = w_n^*(\xi) = n$ for almost all real numbers $\xi$
(with respect to the Lebesgue measure on $\R$),
while $w_n (\xi) = w_n^*(\xi) = {n-1\over 2}$ for almost all
complex numbers (with respect to the Lebesgue measure on $\C$).

Schmidt \cite{SchmAM}
confirmed that with respect to approximation
by algebraic numbers of degree at most $n$, {\it real}
algebraic numbers of degree larger than $n$ behave like almost all real numbers.
Precisely, for every real algebraic number $\xi$ of degree $d$, we have
$$
w_n (\xi) = w_n^* (\xi) = \min\{d-1, n\}  \eqno (1.2)
$$
for every integer $n\geq 1$. The $d-1$ in the right-hand side of
(1.2) is an immediate consequence of the Liouville inequality. A comparison
with Sprind\v zuk's result gives
that if $\xi$ is a real algebraic number of degree $>n$
then $w_n(\xi )=w_n(\eta )$ for almost all $\eta\in\R$, that is,
real algebraic numbers of degree $>n$ are equally well approximable
by algebraic numbers of degree at most $n$ as almost all real numbers.

In this paper we consider the problem to compute $w_n(\xi )$ and $w_n^*(\xi )$
for {\it complex, non-real} algebraic numbers $\xi$.
It follows again from the Liouville inequality that for complex,
non-real algebraic numbers $\xi$
of degree $d\leq n$ one has $w_n(\xi )=w_n^*(\xi )=(d-2)/2$,
but there is no literature about the case where $\xi$ has degree
$d>n$. This case is treated in the present paper.

Our results may be summarized as follows.
Let $\xi$ be a complex, non-real algebraic
number of degree larger than $n$.
Then if $n$ is odd, we have
$w_n(\xi )=w_n^*(\xi )={n-1\over 2}$, while if $n$ is even
we have $w_n(\xi )=w_n^*(\xi )\in \{ {n-1\over 2}, {n\over 2}\}$.
Further, for every even $n$ both cases may occur.
In fact, we are able to decide for every positive even integer $n$
and every complex algebraic number $\xi$ whether
$w_n(\xi )=w_n^*(\xi )= {n-1\over 2}$ or ${n\over 2}$,
except when $n\geq 6$,
$n+2<\deg\xi\leq 2n-2$,
$[\Q (\xi ):\Q (\xi )\cap\R ]=2$, and $1, \xi+\ovxi ,\xi\cdot\ovxi$
are linearly independent over $\Q$.

A comparison
with Sprind\v zuk's result for complex numbers mentioned above
gives that for every even integer $n\geq 2$
there are complex algebraic numbers $\xi$ of degree $>n$
such that $w_n(\xi )> w_n(\eta )$
for almost all complex numbers $\eta$.
So an important consequence of our results is that in contrast to the real case,
for every even integer $n\geq 2$
there are complex algebraic numbers $\xi$ of degree larger than $n$
that are better approximable by algebraic numbers of degree at most $n$
than almost all complex numbers.

We also study how well complex algebraic numbers
can be approximated by algebraic {\it integers} of bounded degree,
and our results support the expectation that complex algebraic numbers
which are unusually well approximable by algebraic numbers of degree at most
$n$, are unusually badly approximable by algebraic integers of degree
at most $n+1$.

We define quantities $\tw_n(\xi )$, $\tw_n^*(\xi )$ analogously
to $w_n(\xi )$, $w_n^*(\xi )$, except that now the approximation
is with respect to monic polynomials in $\Z [X]$ of degree at most $n+1$
and complex algebraic integers of degree at most $n+1$, instead of
polynomials in $\Z [X]$ of degree at most $n$ and complex algebraic
numbers of degree at most $n$.
We prove that
if $\xi$ is a complex algebraic number
of degree
larger than $n$, then
$\tw_n(\xi )=\tw_n^*(\xi )={n-1\over 2}$ if $w_n(\xi )={n-1\over 2}$,
while $\tw_n(\xi )=\tw_n^*(\xi )={n-2\over 2}$ if $w_n(\xi )={n\over 2}$.

Similarly to the case that the number $\xi$ is real algebraic,
in our proofs we apply Schmidt's Subspace Theorem
and techniques from the geometry of numbers.
In this way, we
reduce our approximation problem to a purely algebraic problem
which does not occur in the real case and which leads to additional
difficulties.

\vskip 6mm

\centerline{\bf 2. Main results}

\vskip 5mm

The exponents $w_n$ and $w^*_n$ defined
in the Introduction measure the quality of algebraic approximation, but
do not give any information regarding the number, or the density, of
very good approximations. This lead the authors of \cite{BuLau} to
introduce exponents of {\it uniform} Diophantine approximation.
For a complex number $\xi$
and an integer $n \ge 1$,
we denote by ${\hat w}_n(\xi)$ the supremum of those
real numbers $w$ for which, for every sufficiently large
integer $H$, the inequality
$$
0 < |P(\xi)| \le H^{-w}
$$
is satisfied by an integer polynomial $P$
of degree at most $n$ and height at most $H$.

Khintchine \cite{Kh} proved that ${\hat w}_1 (\xi) = 1$ 
for all irrational real numbers $\xi$. 
Quite unexpectedly, there are real numbers $\xi$ with ${\hat w}_2 (\xi) > 2$.
This was established very recently by Roy \cite{RoyA,RoyB}
(in fact with $\what_2(\xi )={3+\sqrt{5}\over 2}$).  
However, it is still open whether 
there exist an integer $n \ge 3$ and a real number $\xi$ 
such that ${\hat w}_n (\xi) > n$.

Our results show that the three functions
$w_n$, $w^*_n$ and ${\hat w}_n$ coincide on the set of complex
algebraic numbers.
Our first result is as follows.

\proclaim Theorem 1. Let $n$ be a positive integer,
and $\xi$ a complex, non-real algebraic number of degree $d$. Then
$$
\eqalignno{
&w_n(\xi)=w_n^*(\xi )=\what_n(\xi )={d-2\over 2}\ \ \hbox{if $d\leq n+1$,}
&(2.1)\cr
&w_n(\xi)=w_n^*(\xi )=\what_n(\xi )={n-1\over 2}\ \
\hbox{if $d\geq n+2$ and $n$ is odd,}&(2.2)\cr
&w_n(\xi)=w_n^*(\xi )=\what_n(\xi )\in\Big\{{n-1\over 2},{n\over 2}\Big\}\ \
\hbox{if $d\geq n+2$ and $n$ is even.}&(2.3)\cr
}
$$

Thus, Theorem 1 settles completely the case when $n$ is odd.
Henceforth we assume that $n$ is even. In Theorem 2 we give some cases
where $w_n(\xi )={n/2}$ and in Theorem 3 some cases where
$w_n(\xi )={n-1\over 2}$. Unfortunately, we have not been able to compute
$w_n(\xi )$ in all cases. We denote by $\overline{\alpha}$ the complex
conjugate of a complex number $\alpha$.

\proclaim Theorem 2.
Let $n$ be an even positive integer and $\xi$ a
complex, non-real algebraic number of degree $\geq n+2$.
Then $w_n(\xi )=w_n^*(\xi )=\what_n(\xi )={n\over 2}$
in each of the following two cases:
\vskip2pt\noindent
{\bf (i).} $1$, $\xi +\ovxi$ and $\xi\cdot\ovxi$
are linearly dependent over $\Q$;
\vskip2pt\noindent
{\bf (ii).} $\deg\xi =n+2$ and $[\Q (\xi ):\Q (\xi )\cap\R ]=2$.

One particular special case of (i) is when
$\xi = \sqrt{- \alpha}$ for
some positive real algebraic number $\alpha$ of degree $\geq {n\over 2}+1$.
Then $\xi +\ovxi =0$ and so
$w_n(\xi) = w^*_n (\xi) = \hw_n (\xi)=n/2$.

We do not know whether Theorem 2 covers all cases where $w_n(\xi )={n\over 2}$.
We now give some cases where $w_n(\xi )={n-1\over 2}$.

\proclaim Theorem 3.
Let again $n$ be an even positive integer and $\xi$ a
complex, non-real algebraic number of degree $\geq n+2$.
Then $w_n(\xi )=w_n^*(\xi )=\what_n(\xi )={n-1\over 2}$
in each of the following two cases:
\vskip2pt\noindent
{\bf (i).} $[\Q (\xi ):\Q(\xi )\cap\R ]\geq 3$;
\vskip2pt\noindent
{\bf (ii).} $\deg\xi >2n-2$ and $1$, $\xi +\ovxi$, $\xi\cdot\ovxi$
are linearly independent over $\Q$.

For $n=2,4$ we have $2n-2\leq n+2$, so in that case Theorems 2 and 3 cover
all complex algebraic numbers $\xi$. 
Further, for $n=2$, case (ii) of Theorem 2 is implied by case (i).
This leads to the following corollary.

\proclaim Corollary 1. Let $\xi$ be a complex, non-real algebraic number.
\vskip2pt\noindent
{\bf (i).} If $\xi$ has degree $>2$, then
$$
\matrix{
w_2(\xi )=w_2^*(\xi )=\hw_2(\xi )=1\hfill&
\hbox{if $1,\xi +\ovxi, \xi\cdot\ovxi$
are linearly dependent over $\Q$,}\hfill\cr
&\quad\hfill\cr
w_2(\xi )=w_2^*(\xi )=\hw_2(\xi )={1\over 2}\hfill\,\, &
\hbox{otherwise.}\hfill
}
$$
\vskip1pt\noindent
{\bf (ii).} If $\xi$ has degree $>4$, then
$$
\matrix{
w_4(\xi )=w_4^*(\xi )=\hw_4(\xi )=2\hfill&
\hbox{if $1,\xi +\ovxi, \xi\cdot\ovxi$ are linearly dependent over $\Q$}
\hfill\cr
&\hbox{or if $\deg\xi =6$ and $[\Q (\xi ):\Q (\xi )\cap\R ]=2$,}\hfill\cr
&\quad\hfill\cr
w_4(\xi )=w_4^*(\xi )=\hw_4(\xi )={3\over 2}\,\, & \hbox{otherwise.}\hfill
}
$$

Theorems 1,2,3 and Corollary 1 allow us to determine
$w_n(\xi )$, $w_n^*(\xi )$, $\what_n(\xi )$ for every positive integer $n$
and every complex, non-real algebraic number $\xi$, with the exception
of the following case:
\vskip2pt\noindent
{\it $n$ is an even integer with $n\geq 6$,\hfill\break
$\xi$ is a complex algebraic number
such that $n+2<\deg\xi\leq 2n-2$, $[\Q (\xi ):\Q (\xi )\cap\R ]=2$
and $1,\xi +\ovxi ,\xi\cdot\ovxi$ are linearly independent over $\Q$.}
\vskip4mm

We deduce Theorems 1,2,3 from Theorem 4 below.
To state the latter,
we have to introduce some notation.
For $n\in\Z_{>0}$, $\xi\in\C^*$, $\mu\in\C^*$, define the $\Q$-vector space
$$
V_n(\mu ,\xi ):= \{f \in \Q[\variable ] :\, \deg f \leq n,\, \mu f(\xi)\in\R\},
\eqno (2.4)
$$
and for $n\in\Z_{>0}$, $\xi\in\C^*$ denote by $t_n(\xi )$
the maximum over $\mu$ of the dimensions of these spaces, i.e.,
$$
t_n (\xi) := \max\{ \dim_{\Q} V_n(\mu ,\xi ):\, \mu \in \C^* \}.
\eqno (2.5)
$$
It is clear that $t_n(\xi )\leq n+1$
and $t_n(\xi )=n+1$ if and only if $\xi\in\R$.

\proclaim Theorem 4.
Let $n$ be a positive integer and $\xi$ a complex, non-real
algebraic number of degree $>n$. Then
$$
w_n(\xi) = w^*_n (\xi) = {\hat w}_n (\xi) =
\max \biggl\{ {n-1  \over 2}, t_n (\xi) - 1 \biggr\}.
$$

The proof of Theorem 4 is based on Schmidt's Subspace Theorem
and geometry of numbers.
It should be noted that Theorem 4 reduces the problem to determine
how well $\xi$ can be approximated
by algebraic numbers of degree at most $n$ to the algebraic problem
to compute $t_n(\xi )$.
We deduce Theorems 1,2 and 3
by combining Theorem 4 with some properties of the quantity $t_n(\xi )$
proved below.

\vskip 6mm

\centerline{\bf 3. Approximation by algebraic integers}

\vskip 5mm

In view of a transference lemma relating
uniform homogeneous approximation to inhomogeneous
approximation (see \cite{BuLauB}),
for any integer $n \ge 2$, the real numbers $\xi$
with $\hw_n (\xi) > n$ are good candidates for being
unexpectedly badly approximable by algebraic integers of
degree less than or equal to $n+1$.
This has been
confirmed by Roy \cite{RoyC} for the case $n=2$.
Namely, in \cite{RoyB} he proved that there exist real numbers $\xi$ with 
$\what_2(\xi )={3+\sqrt{5}\over 2}>2$, 
and in \cite{RoyC} he used this to prove
that there exist real numbers $\xi$
with the property that $|\xi -\alpha |\gg H(\alpha )^{-(3+\sqrt{5})/2}$
for every algebraic integer $\alpha$ of degree at most $3$.
By a result of Davenport and Schmidt \cite{DaSc69}, the exponent
${3+\sqrt{5}}\over 2$ is optimal.
On the other hand Bugeaud and Teuli\'e \cite{BuTe} proved that
for every $\kappa <3$ and almost all $\xi\in\R$, the inequality 
$|\xi -\alpha |< H(\alpha )^{-\kappa}$ has 
infinitely many solutions in algebraic integers of degree $3$. 

Analogously to the real case one should expect
that complex numbers $\xi$ with $\what_n(\xi )>{n-1\over 2}$
are unusually badly approximable by algebraic integers of degree at most $n+1$.
In Theorem 5 below we confirm this for complex algebraic numbers.
 
We introduce the following quantities
for complex numbers $\xi$ and integers $n\ge 1$:
\vskip0.2cm\noindent
$\bullet$ $\tw_n(\xi )$ denotes the supremum of those real numbers $\tw$ such that
$$
0 < |P(\xi)| \le H(P)^{-\tw}
$$
is satisfied by infinitely many monic polynomials $P\in\Z [\variable ]$
of degree at most $n+1$;
\vskip 0.1cm\noindent
$\bullet$ $\tw_n^*(\xi )$ denotes the supremum of those real numbers $\tw^*$ for which
$$
0 < |\xi -\alpha |\leq H(\alpha )^{-\tw^*-1}
$$
holds for infinitely
many algebraic integers of degree at most $n+1$;
\vskip 0.1cm\noindent
$\bullet$ ${\hat \tw}_n (\xi)$ denotes the supremum of those real numbers $\tw$
with the property that for every sufficiently large
real $H$, there exists a monic integer polynomial $P$ of degree
at most $n+1$ and height at most $H$ such that
$$
0 < |P(\xi)| \le H^{-\tw}.
$$

It is known that every real algebraic number $\xi$ of degree $d$ satisfies
$$
\tw_n (\xi) = \tw_n^* (\xi) = {\hat \tw}_n (\xi) = \min\{d-1, n\}
$$
for every integer $n$ (see \cite{SchmEM,BuLiv}).
Furthermore, methods developed by Bugeaud and Teuli\'e \cite{BuTe} and
Roy and Waldschmidt \cite{RoyWa} allow one to show that for every
positive integer $n$ we have
$$
\eqalign{
&\tw_n (\xi) = \tw_n^* (\xi) = {\hat \tw}_n (\xi) = n
\ \hbox{for almost all $\xi\in\R$},\cr
&\tw_n (\xi) = \tw_n^* (\xi) = {\hat \tw}_n (\xi) = {n-1\over 2}
\ \hbox{for almost all $\xi\in\C$}.\cr
}
$$

We show that for every positive integer $n$
the functions $\tw_n$, $\tw_n^*$,
${\hat\tw}_n$ coincide on the  complex algebraic numbers
and, moreover, that a complex algebraic number $\xi$ is unusually
badly approximable by algebraic integers of degree at most $n+1$
(i.e., has $\tw_n(\xi )=\tw_n^*(\xi )={\hat \tw}_n(\xi )<{n-1\over 2}$)
if and only if it is unusually well approximable by algebraic numbers
of degree at most $n$
(i.e., has $w_n(\xi )=w_n^*(\xi )=\hw_n(\xi )>{n-1\over 2}$).
More precisely,
we prove the following.

\proclaim Theorem 5.
Let $n$ be a positive integer and $\xi$ a complex, non-real algebraic
number of degree $d$.
Then
$$
\eqalignno{
&\tw_n(\xi) = \tw^*_n (\xi) = {\hat \tw}_n (\xi) =
{d-2\over 2}\ \
\hbox{if $d\leq n+1$,}&(3.1)\cr
&\tw_n(\xi )=\tw_n^*(\xi )={\hat \tw}_n(\xi )={n-1\over 2}\
\hbox{if $d\geq n+2$ and $n$ is odd,}&(3.2)\cr
&\tw_n(\xi )=\tw_n^*(\xi )={\hat \tw}_n(\xi )\in
\left\{ {n-2\over 2}, {n-1\over 2}\right\}
\hbox{if $d\geq n+2$ and $n$ is even.}&(3.3)\cr
}
$$
Moreover, if $d\geq n+2$ and $n$ is even then
$$
\tw_n(\xi )=\tw_n^*(\xi )={\hat \tw}_n(\xi )={n-2\over 2}\Longleftrightarrow
w_n(\xi )=w_n^*(\xi )=\hw_n(\xi )={n\over 2}.
$$

Combining Theorem 5 with Corollary 1, we get at once 
the following statement. 

\proclaim Corollary 2. Let $\xi$ be a complex, non-real algebraic number.
\vskip2pt\noindent
{\bf (i).} If $\xi$ has degree $>2$, then
$$
\matrix{
\tw_2(\xi )=\tw_2^*(\xi )={\hat \tw}_2(\xi )=0\hfill&
\hbox{if $1,\xi +\ovxi, \xi\cdot\ovxi$
are linearly dependent over $\Q$,}\hfill\cr
&\quad\hfill\cr
\tw_2(\xi )=\tw_2^*(\xi )={\hat \tw}_2(\xi )={1\over 2}\hfill\,\, &
\hbox{otherwise.}\hfill
}
$$
\vskip1pt\noindent
{\bf (ii).} If $\xi$ has degree $>4$, then
$$
\matrix{
\tw_4(\xi )=\tw_4^*(\xi )={\hat \tw}_4(\xi )=1\hfill&
\hbox{if $1,\xi +\ovxi, \xi\cdot\ovxi$ are linearly dependent over $\Q$}
\hfill\cr
&\hbox{or if $\deg\xi =6$ and $[\Q (\xi ):\Q (\xi )\cap\R ]=2$,}\hfill\cr
&\quad\hfill\cr
\tw_4(\xi )=\tw_4^*(\xi )={\hat \tw}_4(\xi )={3\over 2}\,\, & \hbox{otherwise.}\hfill
}
$$

\vskip 6mm

\goodbreak

\centerline{{\bf 4. Deduction of Theorem 1 from Theorem 4}}

\vskip 5mm

For every positive integer $m$ we define the $\Q$-vector space
$$
W_m:=\{ f\in\Q [\variable ]:\, \deg f\leq m\}
$$
and for any subset $S$ of the polynomial ring $\Q [\variable ]$
and any polynomial $g\in \Q[\variable ]$, we define the set
$g\cdot S :=\{ gf:\, f\in S\}$.

In this section, $n$ is a positive integer, and $\xi$ a complex,
non-real algebraic number of degree $d>n$.
We prove some lemmata about the quantity $t_n(\xi )$ which
in combination with Theorem 4 will imply Theorem 1.
Choose $\mu_0\in\C^*$ such that $\dim V_n(\mu_0,\xi )=t_n(\xi )$.

\proclaim Lemma 4.1.
Let $\mu\in\C^*$ be such that $\dim V_n(\mu ,\xi )>{n+1\over 2}$.
Then $V_n(\mu ,\xi )=V_n(\mu_0,\xi )$.

\proof
Our assumption on $\mu$ clearly implies that $t_n(\xi )>{n+1\over 2}$.
Both vector spaces $V_n(\mu ,\xi )$, $V_n(\mu_0,\xi )$ are contained
in the same $n+1$-dimensional vector space, hence they have non-zero
intersection. Let $f_1\in\Q [\variable ]$ be a non-zero polynomial lying in both spaces
and put $\mu_1 := f_1(\xi )^{-1}$. Then $\mu_1/\mu\in\R$,
$\mu_1/\mu_0\in\R$, hence $V_n(\mu ,\xi )=V_n(\mu_1,\xi )=V_n(\mu_0,\xi )$.
\qed

\proclaim Lemma 4.2.
Suppose that $t_n(\xi )>{n+1\over 2}$. Then
\vskip2pt\noindent
{\bf (i).} $W_{n+1}$ is the direct sum of the $\Q$-vector spaces
$V_n(\mu_0,\xi )$ and $X\cdot V_n(\mu_0,\xi)$.
\vskip2pt\noindent
{\bf (ii).} $n$ is even, $t_n(\xi )={n+2\over 2}$.

\proof
Suppose that $V_n(\mu_0,\xi )\cap X\cdot V_n(\mu_0,\xi )\not=\{ 0\}$.
Choose a non-zero polynomial $f$ in the intersection of both spaces.
Then $f=Xg$ where $g\in V_n(\mu_0,\xi )$. Hence
$$
\xi ={f(\xi )\over g(\xi )}={\mu_0 f(\xi )\over \mu_0 g(\xi )}\in\R \, ,
$$
which is against our assumption. Therefore,
$V_n(\mu_0,\xi )\cap X\cdot V_n(\mu_0,\xi )=\{ 0\}$.
From our assumption on $\xi$ it follows that $t_n(\xi )\geq {n+2\over 2}$.
Further, both $V_n(\mu_0,\xi )$
and $X\cdot V_n(\mu_0,\xi )$ are linear subspaces of $W_{n+1}$.
Hence by comparing dimensions,
$$
2\cdot{n+2\over 2}\leq 2t_n(\xi )
=\dim \Big(V_n(\mu_0,\xi )+X\cdot V_n(\mu_0,\xi )\Big)\leq \dim W_{n+1} =n+2.
$$
This implies (i) and (ii).
\qed

\proclaim Lemma 4.3.
Let $\xi$ be a complex, non-real algebraic number of degree $d>1$.
Then $t_{d-1}(\xi )\leq {d\over 2}$.

\proof
Choose $\mu_0\in\C^*$ such that $\dim V_{d-1}(\mu_0,\xi )=t_{d-1}(\xi )$.
Pick a non-zero polynomial $f_0\in V_{d-1}(\mu_0,\xi )$. Then for
every $f\in V_{d-1}(\mu_0,\xi )$ we have
${f(\xi )\over f_0(\xi )}
={\mu_0f(\xi )\over \mu_0f_0(\xi )}\in \Q (\xi )\cap\R$.
For linearly independent polynomials $f\in \Q [\variable ]$
of degree at most $d-1=\deg\xi -1$, the corresponding
quantities $f(\xi )/f_0(\xi )$ are linearly independent over $\Q$.
Hence $t_{d-1}(\xi )\leq [\Q (\xi )\cap\R :\Q ]\leq {d\over 2}$.
\qed

In the proof of Theorem 1 we use the following observations.

\proclaim Lemma 4.4.
Let $\xi$ be a complex number and $n$ a positive integer. Then
\hfill\break
{\bf (i).} $w_n^*(\xi )\leq w_n(\xi )$,
\hfill\break
{\bf (ii).} $\hw_n(\xi )\leq w_n(\xi )$.

\proof
If $\alpha$ is an algebraic number of degree $n$ with minimal polynomial
$P\in\Z [\variable ]$, we have 
$|P(\xi )|\ll H(P) \cdot \min\{1,|\alpha -\xi |\}$,
where the implied
constant depends only on $\xi$ and on $n$.
This implies (i).
If for some $w\in\R$ there exists $H_0$ such that for every $H\geq H_0$
there exists an integer polynomial $P$
of degree at most $n$ with
$0<|P(\xi )|\leq H^{-w}$, $H(P)\leq H$, then clearly, there are infinitely
many integer polynomials $P$
of degree at most $n$ such that
$0< |P (\xi )|\leq H(P)^{-w}$. This implies (ii).
\qed

\noindent
{\bf Proof of Theorem 1:}
Constants implied by $\ll$ and $\gg$ depend only on $n$, $\xi$.
We first prove (2.1). Assume that $d\leq n+1$.
In view of Lemma 4.4,
it suffices to prove that
$$
w_n(\xi )\leq {d-2\over 2},\ \
w_n^*(\xi )\geq  {d-2\over 2},\ \  \hw_n(\xi )\geq {d-2\over 2}.
$$
To prove the former, denote by $\xi^{(1)},\ldots ,\xi^{(d)}$
the conjugates of $\xi$, where $\xi^{(1)}=\xi$, $\xi^{(2)}=\overline{\xi}$.
For some $a\in\Z_{>0}$, the polynomial
$Q:= a\prod_{i=1}^d (X-\xi^{(i)})$ has integer coefficients,
and for any polynomial $P\in\Z [\variable ]$ of degree at most $n$
with $P(\xi )\not= 0$,
the resultant $R(P,Q)=a^n\prod_{i=1}^d P(\xi^{(i)})$ is a non-zero
rational integer.
This gives the Liouville inequality
$$
|P(\xi )|^2=|P(\xi )P(\overline{\xi})|
\gg { |R(P,Q)|\over |P(\xi^{(3)})\cdots P(\xi^{(d)})|}
\gg H(P)^{2-d}.\eqno (4.1)
$$
Consequently, $w_n(\xi )\leq {d-2\over 2}$.

By Theorem 4 with $n=d-1$ and by Lemma 4.3 we have
$w_{d-1}^*(\xi )=\hw_{d-1}(\xi )={d-2\over 2}$.
Hence for $n\geq d$ we have
$$
w_n^*(\xi )\geq w_{d-1}^*(\xi )= {d-2\over 2},\ \
\hw_n(\xi )\geq \hw_{d-1}(\xi )= {d-2\over 2}.
$$
This completes the proof of (2.1).

Statements (2.2), (2.3)
follow immediately by combining Theorem 4 with part (ii) of Lemma 4.2.
This completes the proof of Theorem 1.
\qed

\vskip 6mm

\centerline{{\bf 5. Deduction of Theorem 2 from Theorem 4}}

\vskip 5mm

To deduce Theorem 2 from Theorem 4,
we prove again the necessary properties
for the quantity $t_n(\xi )$ defined by (2.5).

\proclaim Lemma 5.1.
Assume that $n$ is even, and that $\xi$ is a complex, non-real
algebraic number of degree $>n$ such that
$1$, $\xi +\ovxi$ and $\xi\cdot \ovxi$
are linearly dependent over $\Q$. Then
$$
t_n(\xi )= {n+2\over 2}.
$$

\proof
We use the easy observation that $t_n(\xi +c)=t_n(\xi )$ for any $c\in\Q$.

Put $\beta := \xi +\ovxi$, $\gamma := \xi\cdot\ovxi$.
Our assumption on $\xi$ implies that either $\beta\in\Q$, or $\gamma =a+b\beta$
for some $a,b\in\Q$. By our observation,
the first case can be reduced to
$\beta =0$ by replacing $\xi$ by $\xi -{1\over 2}\beta$.
Then $\xi =\sqrt{-\gamma}$ with $\gamma >0$. Likewise,
the second case can be reduced to $\gamma =a\in\Q$ by replacing $\xi$
by $\xi -b$. Then $\xi ={1\over 2}\Big(\beta \pm\sqrt{\beta^2 -4a}\Big)$
with $a\in\Q$ and $a>\beta^2/4$.
\vskip2pt

{\bf Case I.} $\xi =\sqrt{-\gamma}$ with $\gamma >0$.
\vskip 2pt\noindent
In this case,
$$
V_n(1,\xi )=\{ f\in\Q [\variable ]:\, \deg f\leq n,\, f(\xi )\in\R\}
=\biggl\{ \sum_{i=0}^{n/2} c_iX^{2i}:\, c_0,\ldots ,c_{n/2}\in\Q\biggr\}.
$$
So $t_n(\xi )\geq \dim V_n(1 ,\xi )={n+2\over 2}$.
Hence by Lemma 4.2 we have
$t_n(\xi )={n+2\over 2}$.
\vskip2pt

{\bf Case II.} $\gamma =\xi\cdot\ovxi =a\in\Q^*$.
\vskip 2pt\noindent
Put $\mu :=\xi^{-n/2}$. Then for a polynomial $f=\sum_{i=0}^n c_iX^i\in \Q [\variable ]$
we have, recalling our assumption that $\xi$ has degree larger than $n$,
$$
\eqalign{
\mu f(\xi )\in\R &\Longleftrightarrow \xi^{-n/2}f(\xi )=\ovxi^{(-n/2)}f(\ovxi )
\Longleftrightarrow \xi^{-n/2}f(\xi ) =(a/\xi )^{-n/2}f(a/\xi )\cr
&\Longleftrightarrow a^{n/2}f(\xi )= \xi^nf(a/\xi )\Longleftrightarrow
a^{n/2}f(X)=X^nf(a/X)\cr
&\Longleftrightarrow a^{n/2}c_i =a^{n-i}c_{n-i}\
\hbox{for $i=0,\ldots ,n$.}
}
$$
This implies $t_n(\xi )\geq \dim V_n(\mu ,\xi )={n+2\over 2}$.
Hence $t_n(\xi )={n+2\over 2}$ in view of Lemma 4.2.
\qed

\proclaim Lemma 5.2.
Let $n$ be an even positive integer,
and $\xi$ a complex algebraic number of degree $n+2$.
Suppose that $\, [\Q (\xi ):\Q (\xi )\cap\R ]=2$. Then
$$
t_n(\xi )={n+2\over 2}.
$$

\proof
Write $k:= n/2$.
Then $\Q (\xi )\cap\R$ has degree $k+1$. We prove that there exists
$\mu\in \Q (\xi )^*$ such that $\dim V_n(\mu ,\xi )\geq k+1 ={n+2\over 2}$.
Then from Lemma 4.2 it follows that $t_n(\xi )={n+2\over 2}$.

Let $\{\omega_1,\ldots ,\omega_{k+1}\}$ be a $\Q$-basis of $\Q (\xi )\cap\R$.
Then $\omega_1,\ldots ,\omega_{k+1}$, $\xi\omega_1,\ldots ,\xi\omega_{k+1}$
form a $\Q$-basis of $\Q (\xi )$, every element of $\Q (\xi )$ can be expressed
uniquely as a $\Q$-linear combination of these numbers,
and a number in $\Q (\xi )$ thus expressed belongs to $\Q (\xi )\cap\R$
if and only if its coefficients with respect to
$\xi\omega_1,\ldots ,\xi\omega_{k+1}$ are $0$.

For $i,j=0,\ldots ,2k+1$ we have
$$
\xi^{i+j}=\sum_{l=1}^{k+1} a_{ij}^{(l)}\omega_l +
\sum_{l=1}^{k+1} b_{ij}^{(l)}\xi\omega_l\quad
\hbox{with $a_{ij}^{(l)},b_{ij}^{(l)}\in\Q$.}
$$
Write $\mu\in \Q (\xi )$ as $\mu =\sum_{i=0}^{2k+1} u_i\xi^i$
with $u_0,\ldots ,u_{2k+1}\in\Q$
and write $f\in V_n(\mu ,\xi )$ as $f=\sum_{j=0}^{2k} x_jX^j$
with $x_0,\ldots ,x_{2k}\in\Q$. 
Then
$$
\mu f(\xi )=\sum_{l=1}^{k+1} \omega_l
\left\{ \sum_{j=0}^{2k}\Big( \sum_{i=0}^{2k+1} a_{ij}^{(l)}u_i\Big)x_j\right\}
+
\sum_{l=1}^{k+1} \xi\omega_l\left\{ \sum_{j=0}^{2k}
\Big( \sum_{i=0}^{2k+1} b_{ij}^{(l)}u_i\Big)x_j\right\}.
$$
So $f=\sum_{j=0}^{2k} x_jX^j\in V_n(\mu ,\xi )$, i.e.,
$\mu f(\xi )\in \Q (\xi )\cap\R$, if and only if
$$
L_{\mu}^{(l)}({\bf x}):=\sum_{j=0}^{2k}
\Big( \sum_{i=0}^{2k+1} b_{ij}^{(l)}u_i\Big)x_j\,
=\, 0\quad \hbox{for $l=1,\ldots ,k+1$,}
\eqno (5.1)
$$
where ${\bf x}=(x_0,\ldots , x_{2k})$.

We choose $\mu\in\Q (\xi )^*$ to make one of the linear forms in (5.1),
for instance $L_{\mu}^{(k+1)}$, vanish identically. This amounts to
choosing a non-zero vector ${\bf u}=(u_0,\ldots ,u_{2k+1})\in\Q^{2k+2}$
such that
$$
\sum_{i=0}^{2k+1} b_{ij}^{(k+1)}u_i =0\quad\hbox{for $j=0,\ldots ,2k$.}
$$
This is possible since a system of $2k+1$ linear equations in $2k+2$
unknowns has a non-trivial solution. Thus, (5.1) becomes a system of
$k$ equations in $2k+1$ unknowns over $\Q$, and the solution space of this
system has dimension at least $k+1$.
Consequently, $V_n(\mu ,\xi )$ has dimension at least $k+1={n+2\over 2}$.
This proves Lemma 5.2.
\qed

Now Theorem 2 follows at once by combining Theorem 4 with Lemmata 5.1 and 5.2.
\qed

\vskip 6mm

\centerline{{\bf 6. Deduction of Theorem 3 from Theorem 4}}

\vskip 5mm

We prove some results about the quantity $t_n(\xi )$
which, in combination with
Theorem 4, will yield Theorem 3.

\proclaim Lemma 6.1.
Let $n$ be an even positive integer and
$\xi$ a complex, non-real algebraic number of degree $>n$.
Assume that $t_n(\xi )>{n+1\over 2}$.
\vskip2pt\noindent
{\bf (i).} $[\Q (\xi ):\Q (\xi )\cap\R ]=2$.
\vskip2pt\noindent
{\bf (ii).} If moreover $\deg\xi >2n-2$, then
$1$, $\xi +\ovxi$, $\xi\cdot\ovxi$ are linearly dependent over $\Q$.

\proof
Put $\beta :=\xi +\ovxi$, $\gamma := \xi\cdot\ovxi$.
Choose $\mu_0$ such that $\dim V_n(\mu_0,\xi )=t_n(\xi )$.
By part (i) of Lemma 4.2, every polynomial in $\Q [\variable ]$ of degree at most $n+1$
can be expressed uniquely as a sum of a polynomial in $V_n(\mu_0,\xi )$
and a polynomial in $X\cdot V_n(\mu_0,\xi )$.
In particular, for every
non-zero polynomial $f\in V_n (\mu_0,\xi )$ of degree $\leq n-1$, there are
polynomials $g,h\in V_n(\mu_0,\xi )$, uniquely determined by $f$,
such that
$$
X^2f=Xg+h.
\eqno (6.1)
$$
This implies that $\xi$ is a zero of the polynomial
$X^2-(g(\xi )/f(\xi ))X -(h(\xi )/f(\xi ))$. On the other hand, there is
a unique monic quadratic polynomial with real coefficients having $\xi$
as a zero, namely $X^2-\beta X+\gamma$, and
$$
{g(\xi )\over f(\xi )}={\mu_0g(\xi )\over \mu_0f(\xi )}\in\R,\quad
{h(\xi )\over f(\xi )}={\mu_0h(\xi )\over \mu_0f(\xi )}\in\R.
$$
Therefore,
$$
{g(\xi )\over f(\xi )}=\beta,\quad {h(\xi )\over f(\xi )}=-\gamma\, .
\eqno (6.2)
$$
So $\beta ,\gamma\in\Q(\xi )\cap\R$. This implies (i).

To prove (ii), we proceed by induction on $n$.
First let $n=2$.
By assumption, there is $\mu_0\in\C^*$ such that $V_2(\mu_0 ,\xi )$
has dimension larger than $1$. This means that
there are non-zero polynomials $f_1,f_2\in V_2(\mu_0 ,\xi )$
with $\deg f_1<\deg f_2\leq 2$.
We have
$f_1(\xi )f_2(\ovxi )\in (\mu_0\overline{\mu_0})^{-1}\R =\R$, hence
$$
f_1(\xi )f_2(\ovxi )-f_1(\ovxi )f_2(\xi )=0.
$$
First suppose that $f_2$ has degree $1$.
Then $f_1$ has degree $0$, therefore,
$f_1=c_1$, $f_2=c_2+c_3X$ with $c_1c_3\not= 0$. Hence
$$
0=f_1(\xi )f_2(\ovxi )-f_1(\ovxi )f_2(\xi )
= c_1c_3(\ovxi -\xi )
$$
which is impossible since $\xi\not\in\R$.
Now suppose that $f_2$ has degree $2$. Then
$f_1=c_1+c_2X$,
$f_2=c_3+c_4X+c_5X^2$ with $c_1,\ldots ,c_5\in\Q$,
hence
$$
0=f_1(\xi )f_2(\ovxi )-f_1(\ovxi )f_2(\xi )
= (\ovxi -\xi )(c_1c_4-c_2c_3 +c_1c_5\beta +c_2c_5\gamma ).
$$
We have $(c_1,c_2)\not= (0,0)$ since $f_1\not= 0$, while
$c_5\not =0$ since $f_2$ has degree $2$, and further
$\xi\not\in\R$. Hence
$1, \beta ,\gamma$ are $\Q$-linearly dependent.

Now let $n$ be an even integer with $n\geq 4$. Assume part (ii)
of Lemma 6.1 is true
if $n$ is replaced by any positive even integer smaller than $n$.
There is $\mu_0\in\C^*$ such that $\dim V_n(\mu_0 ,\xi )=: t>{n+1\over 2}$.
Let $f_1,\ldots ,f_t$ be a basis of $V_n(\mu_0 ,\xi )$ with
$\deg f_1 <\deg f_2<\cdots <\deg f_t\leq n$. So in particular,
$\deg f_{t-1}\leq n-1$.

First assume that $a:={\rm gcd}(f_1,\ldots ,f_{t-1})$ is a polynomial
of degree at least $1$. Let $\tilde{f}_i :=f_i/a$ for $i=1,\ldots ,t-1$.
Put $\tilde{\mu_0}:= \mu_0 a(\xi )$. Then $\tilde{f}_1,\ldots ,\tilde{f}_{t-1}$
are linearly independent polynomials of degree at most $n-2$ with
$\tilde{\mu_0}\tilde{f}_i(\xi )\in\R$ for $i=1,\ldots ,t-1$.
Hence
$$
t_{n-2}(\xi )\geq \dim V_{n-2}(\tilde{\mu_0},\xi )\geq t-1>{(n-2)+1\over 2}.
$$
So by the induction hypothesis, $1,\beta ,\gamma$ are linearly dependent
over $\Q$.

Now assume that ${\rm gcd}(f_1,\ldots ,f_{t-1})=1$.
By (6.1),
for $i=1,\ldots ,t-1$
there are polynomials $g_i,h_i\in V_n(\mu_0 ,\xi )$ such that
$X^2f_i =Xg_i +h_i$ for $i=1,\ldots ,t-1$ and by (6.2) we have
$$
{g_i(\xi )\over f_i(\xi )}=\beta ,\ \ {h_i(\xi )\over f_i(\xi )}=-\gamma
\quad\hbox{for $i=1,\ldots ,t-1$.}
$$
The polynomials $h_i$ are all divisible by $X$. Therefore, $\xi$ is a
common zero of the polynomials
$$
f_i\cdot{h_j\over X}\,-\, f_j\cdot {h_i\over X}\quad (1\leq i,j\leq t-1).
$$
Each of these polynomials has degree at most $2n-2$ and, by assumption,
$\xi$ has degree $>2n-2$. Therefore, these polynomials are all identically $0$.
Since by assumption ${\rm gcd}(f_1,\ldots ,f_{t-1})=1$, this implies that
there is a polynomial $a\in\Q [\variable ]$ with
$h_i/X =af_i$ for $i=1,\ldots , t-1$.

Now $a$ cannot be equal to $0$ since otherwise $\gamma =\xi\cdot\ovxi$
would be $0$ which is impossible. Further, $a$ cannot be a
constant $c\in\Q^*$ since otherwise, we would have
$\xi =h_i(\xi )/cf_i(\xi )=-\gamma /c\in\R$ which is impossible.
Hence $a$ has degree at least $1$.
But then $\deg f_i \leq \deg h_i-2\leq n-2$ for $i=1,\ldots ,t-1$.
This implies
$$
t_{n-2}(\xi )\geq \dim V_{n-2}(\mu ,\xi )\geq t-1>{n-2+1\over 2}.
$$
Now again the induction hypothesis can be applied, and we infer that
$1,\beta ,\gamma$ are linearly dependent over $\Q$.
This completes our proof.
\qed

Theorem 3 follows at once by combining Theorem 4 with Lemma 6.1.
\qed

\vskip 6mm

\centerline{\bf 7. Consequences of the Parametric Subspace Theorem}

\vskip 5mm

In this section we have collected some applications of the Parametric
Subspace Theorem which are needed in both the proofs
of Theorem 4 and Theorem 5. Our arguments are a routine extension of
Chapter VI, \S\S 1,2 of Schmidt's Lecture Notes \cite{SchmLN},
but for lack of a convenient reference we have included the proofs.

We start with some notation.
For a linear form $L=\sum_{i=1}^n \alpha_iX_i$ with complex coefficients,
we write $\Re (L):= \sum_{i=1}^n (\Re\alpha_i)X_i$
and $\Im (L):=\sum_{i=1}^n (\Im\alpha_i)X_i$.
For a linear subspace $U$ of $\Q^n$, we denote by $\R U$ the
$\R$-linear subspace of $\R^n$ generated by $U$.
We say that linear forms
$L_1,\ldots ,L_s$ in $X_1,\ldots ,X_n$ with complex
coefficients are linearly dependent on a linear subspace $U$ of $\Q^n$
if there are complex
numbers $a_1,\ldots ,a_s$, not all zero, such that $a_1L_1+\cdots +a_sL_s$
vanishes identically on $U$. Otherwise, $L_1,\ldots ,L_s$ are said to be
linearly independent on $U$.
\vskip 2pt

Our main tool is the so-called Parametric Subspace Theorem which is stated
in Proposition 7.1 below.
We consider symmetric convex bodies
$$
\Pi (H) :=
\{ {\bf x}\in \R^n :\, |L_i({\bf x})|\leq H^{-c_i}\ (i=1,\ldots ,r)\}
\eqno (7.1)
$$
where $r\geq n$, $L_1,\ldots ,L_r$
are linear forms with real algebraic coefficients
in the $n$ variables $X_1,\ldots ,X_n$, $c_1,\ldots,c_r$ are reals,
and $H$ is a real $\geq 1$. We will refer to $c_i$ as the {\it $H$-exponent
corresponding to $L_i$}.

\proclaim Proposition 7.1.
Assume that there are indices $i_1,\ldots ,i_n\in\{ 1,\ldots ,r\}$ such that
$$
{\rm rank} (L_{i_1},\ldots ,L_{i_n})=n,\ c_{i_1}+\cdots +c_{i_n}>0. \eqno (7.2)
$$
Then there is a finite collection of proper linear subspaces
$\{ T_1,\ldots ,T_t\}$ of $\Q^n$ such that for every $H\geq 1$
there is $T_i\in \{ T_1,\ldots ,T_t\}$ with
$$
\Pi (H)\cap\Z^n\subset T_i\, .
$$

\proof
This is a special case of Theorem 1.1 of \cite{EvSchl},
where a quantitative version was given with an explicit upper bound
for the number of subspaces $t$.
In fact, in its qualitative form
this result was already proved implicitly by Schmidt.
\qed

\proclaim Lemma 7.2.
Let $L_1,\ldots ,L_r$ be linear forms in $X_1,\ldots ,X_n$
with real algebraic coefficients and with
${\rm rank}(L_1,\ldots ,L_r)=n$, let $c_1,\ldots ,c_r$ be reals,
and let $\{ M_1,\ldots ,M_s\}$ be a (possibly empty) collection of
linear forms in $X_1,\ldots ,X_n$
with complex coefficients.\hfill\break
Assume that for every
non-zero linear subspace $U$ of $\Q^n$
on which none of $M_1,\ldots ,M_s$ vanishes identically
there are indices $i_1,\ldots ,i_m\in \{ 1,\ldots ,r\}$ ($m=\dim U$)
such that
$$
\hbox{$L_{i_1},\ldots ,L_{i_m}$ are linearly independent on $U$,}\quad
c_{i_1}+\cdots +c_{i_m}>0.
\eqno (7.3)
$$
Then there is $H_0>1$ such that if
there is ${\bf x}$ with
$$
{\bf x}\in \Pi (H)\cap\Z^n,\ {\bf x}\not= 0,\
M_j({\bf x})\not= 0\ \hbox{for $j=1,\ldots ,s$},
$$
then $H\leq H_0$.

\proof
A subspace $U$ of $\Q^n$ is called admissible if none of
$M_1,\ldots ,M_s$ vanishes identically on $U$.
From (7.3) and Proposition 7.1
it follows easily, that for every
admissible linear subspace $U$ of $\Q^n$ of dimension $\geq 2$,
there is a finite collection
$\{ U_1,\ldots ,U_u\}$ of admissible proper linear subspaces of $U$,
such that for every $H\geq 1$ there is $U_i\in\{ U_1,\ldots ,U_u\}$
with
$$
\{ {\bf x}\in \Pi (H)\cap\Z^n\cap U :\,
M_j({\bf x})\not={\bf 0}\ \hbox{for $j=1,\ldots ,s$}\}\subset U_i\, .
$$

By repeatedly applying this, it follows that there is a finite collection
$\{ V_1,\ldots ,V_v\}$ of admissible one-dimensional linear subspaces
of $\Q^n$, such that for every $H\geq 1$
there is $V_i\in\{ V_1,\ldots ,V_v\}$
with
$$
\{ {\bf x}\in \Pi (H)\cap\Z^n,\
M_j({\bf x})\not={\bf 0}\ \hbox{for $j=1,\ldots ,s$}\}\subset V_i\, .
$$

Let $V$ be one of these subspaces.
Choose a non-zero vector ${\bf x}_0\in V\cap\Z^n$ whose coefficients
have gcd $1$. Such a vector is up to sign uniquely determined by $V$,
and every vector in $V$ is a scalar multiple of ${\bf x}_0$.
By assumption, there is $i\in\{ 1,\ldots ,r\}$ such that
$L_i({\bf x}_0)\not= 0$ and $c_i>0$. Now let ${\bf x}$ be a non-zero
vector in $\Pi (H)\cap\Z^n\cap V$. Then ${\bf x}=a{\bf x}_0$ for some
non-zero integer $a$, hence
$$
H^{-c_i}\geq |L_i({\bf x})|\geq |L_i({\bf x}_0)|,
$$
which implies that $H\leq H_V$ for some finite constant $H_V$ depending
only on $V$.

Now Lemma 7.2 is satisfied with $H_0=\max_{i=1,\ldots ,v} H_{V_i}$.
\qed

Denote by $\lambda_1(H),\ldots ,\lambda_n(H)$ the successive minima
of $\Pi (H)$. Recall that $\lambda_i(H)$ is the minimum of all
positive reals $\lambda$ such that $\lambda \Pi (H)$ contains $i$ linearly
independent points from $\Z^n$.

\proclaim Lemma 7.3.
Let $L_1,\ldots ,L_r$ be linear forms in $X_1,\ldots ,X_n$
with real algebraic coefficients and with
${\rm rank}(L_1,\ldots ,L_r)=n$ and let $c_1,\ldots ,c_r$ be reals.
Put
$$
E := {1\over n}\max \{c_{i_1}+\cdots +c_{i_n}\} \eqno (7.4)
$$
where the maximum is taken over all tuples $i_1,\ldots ,i_n$
such that $L_{i_1},\ldots ,L_{i_n}$ are linearly independent.
\vskip 2pt\noindent
{\bf (i).} There is a constant $c>0$ depending only on $n,L_1,\ldots ,L_r$
such that for every $H\geq 1$ we have $\lambda_1(H)\leq cH^{E}$.
\vskip 2pt\noindent
{\bf (ii).} Assume that for every
non-zero linear subspace $U$ of $\Q^n$
there are indices $i_1,\ldots ,i_m\in \{ 1,\ldots ,r\}$ ($m=\dim U$)
such that
$$
\hbox{$L_{i_1},\ldots ,L_{i_m}$ are linearly independent on $U$,}\quad
{1\over m}(c_{i_1}+\cdots +c_{i_m})\geq E .
\eqno (7.5)
$$
Then for every $\vep >0$ there is $H_{\vep}>1$ such that for every
$H>H_{\vep}$ we have
$$
H^{E -\vep}< \lambda_1(H)\leq \cdots \leq \lambda_n(H)< H^{E  +\vep}.
$$

\proof
In what follows, the constants implied by $\ll$ and $\gg$ may depend
on $L_1,\ldots ,L_r$, $c_1,\ldots ,c_r$, $n$, $\vep$,
but are independent of $H$.
Without loss of generality, $L_1,\ldots ,L_n$ are linearly independent
and $c_1\geq\cdots\geq c_r$.

We first prove (i).
Let $\Pi'(H)$ be the set of ${\bf x}\in\R^n$
with $|L_i({\bf x})|\leq H^{-c_i}$
for $i=1,\ldots ,n$ (so with only $n$ instead of $r$ inequalities).
There is a constant $\lambda_0>0$ such that
$\Pi (H)\supseteq \lambda_0\Pi'(H)$
and this implies at once
$$
{\rm Vol}(\Pi (H))\gg {\rm Vol}(\Pi'(H))\gg H^{-(c_1+\cdots +c_n)}=
H^{-nE }.
$$
So by Minkowski's Theorem on successive minima,
$$
\prod_{i=1}^n \lambda_i(H) \ll H^{nE }.
\eqno (7.6)
$$
This implies (i).

We now prove (ii), and assume that for every non-zero linear subspace
$U$ of $\Q^n$ there are indices $i_1,\ldots ,i_m$ with (7.5).
Let $\vep >0$. We first show that for every sufficiently large $H$
we have
$$
\lambda_1(H)> H^{E -\vep /n},
\eqno (7.7)
$$
in other words, that for every sufficiently large $H$ the convex body
$$
H^{E  -\vep /n}\Pi (H)=
\{ {\bf x}\in\R^n : |L_i({\bf x})|\leq H^{E -c_i-\vep /n}\
(i=1,\ldots ,r)\}
$$
does not contain non-zero points ${\bf x}$ in ${\bf Z}^n$.

We apply Lemma 7.2
with $c_i-E +\vep /n$ instead of $c_i$ for $i=1,\ldots ,r$.
From our assumption it follows that
for every non-zero linear subspace $U$ of $\Q^n$
there are indices $i_1,\ldots ,i_m$ ($m=\dim U$) such that
$L_{i_1},\ldots ,L_{i_m}$ are linearly independent on $U$ and
$$
\sum_{j=1}^m ( c_{i_j}-E+\vep /n )=
(\sum_{j=1}^m c_{i_j})-mE +m\vep /n >0.
$$
So condition (7.3) is satisfied, and therefore we have
$H^{E -\vep/n}\Pi (H)\cap\Z^n=\{ {\bf 0}\}$ for every sufficiently
large $H$. This proves (7.7).

Now a combination of (7.7) with (7.6) immediately gives (ii).
\qed

Let $n$ be a positive integer and $\xi$ a complex, non-real algebraic number
of degree larger than $n$.
Define the linear forms
$$
L_1 := \Re\Big(\sum_{i=0}^n \xi^iX_i\Big),\quad
L_2 := \Im\Big(\sum_{i=0}^n \xi^iX_i\Big)
\eqno (7.8)
$$
and the symmetric convex body
$$
\eqalign{
K(\xi ,n,w,H):= \{ {\bf x}\in\R^{n+1}:\,
|L_1({\bf x})|\leq H^{-w},\
|L_2({\bf x})| & \leq H^{-w}, \cr
& |x_0|\leq H,\ldots ,|x_n|\leq H\},\cr}
\eqno (7.9)
$$
where ${\bf x}=(x_0,\ldots ,x_n)$ and $w\in\R$.
We denote by $\lambda_i (\xi ,n,w,H)$ ($i=1,\ldots ,n+1$)
the successive minima of this body.

Recall that $V_n(\mu ,\xi )$ consists of the
polynomials $f\in\Q [\variable ]$ of degree at most $n$ for which
$\mu f(\xi )\in\R$.
We start with a simple lemma.

\proclaim Lemma 7.4.
{\bf (i).} Let $U$ be a non-zero linear subspace of $\Q^{n+1}$.
Then at least one of the linear forms $L_1$, $L_2$ does not vanish
identically on $U$.
\vskip 2pt\noindent
{\bf (ii).} Let $U$ be a linear subspace of $\Q^{n+1}$.
Then $L_1,L_2$ are linearly
dependent on $U$ if and only if there is $\mu\in\C^*$ such that
$$
U\subset \{ {\bf x}\in\Q^{n+1}:\, \sum_{i=0}^n x_iX^i\in V_n(\mu ,\xi )\}.
$$

\proof
(i). If $L_1,L_2$ would both vanish identically on $U$, then so would
$L_1+\sqrt{-1}\cdot L_2 =\sum_{i=0}^n x_i\xi^i$. But this is
impossible since $\xi$ has degree larger than $n$.

(ii). The linear forms
$L_1,L_2$ are linearly dependent on $U$ if and only if
there are $\alpha ,\beta\in\R$ such that $\alpha L_1 +\beta L_2$
is identically zero on $U$. Using
$$
\alpha L_1({\bf x}) +\beta L_2({\bf x})
= {\rm Im}\Big( \mu\sum_{i=0}^n x_i\xi^i\Big)\ \hbox{with }
\mu = \beta +\sqrt{-1}\cdot \alpha ,
$$
one verifies at once that $L_1,L_2$ are linearly dependent on $U$
if and only if for every ${\bf x}\in U$
the polynomial $\sum_{i=0}^n x_iX^i$ belongs to $V_n(\mu ,\xi )$.
\qed

Let $t_n(\xi )$ be the quantity defined by (2.5).
By Lemma 4.2, we have either $t_n(\xi )\leq {n+1\over 2}$ or
$t_n(\xi )={n+2\over 2}$.
In what follows we have to distinguish between these two cases.
In the proofs below, constants implied by $\ll$ and $\gg$
may depend on $\xi$, $n$, $w$, and on an additional parameter $\vep$,
but are independent of $H$.

\proclaim Lemma 7.5.
Assume that $t_n(\xi )\leq (n+1)/2$ and let $w\geq -1$.
\vskip 2pt\noindent
{\bf (i).} There is a constant $c=c(\xi ,n)>0$ such that for every $H\geq 1$
we have $\lambda_1(\xi ,n,w,H)\leq cH^{2w-n+1\over n+1}$.
\vskip 2pt\noindent
{\bf (ii).}
For every $\vep >0$ there
is $H_{1,\vep}>1$ such that for every $H>H_{1,\vep}$ we have
$$
H^{{2w-n+1\over n+1}-\vep}<\lambda_1(\xi ,n,w,H)\leq\cdots
\leq \lambda_{n+1}(\xi ,n,w,H)<
H^{{2w-n+1\over n+1}+\vep}.
\eqno (7.10)
$$

\proof
In the situation being considered here, for the quantity $E$ defined
by (7.4) we have
$E ={2w-n+1\over n+1}$. Thus, part (i) of Lemma 7.5 follows at once
from part (i) of Lemma 7.3.

We deduce part (ii) of Lemma 7.5 from part (ii) of Lemma 7.3.
and to this end we have to verify the conditions of the latter.
First let $U$ be a linear subspace of $\Q^{n+1}$
of dimension $m> t_n(\xi )$. By part (ii) of Lemma 7.4, the linear forms
$L_1,L_2$ are linearly independent on $U$.
Pick $m-2$
linear forms from $X_0,\ldots ,X_n$ which together with
$L_1,L_2$ are linearly independent on $U$. Then the sum of the $H$-exponents
corresponding to these linear forms is equal to $2w-m+2$, and
$$
{2w -m+2\over m}\geq {2w-n+1\over n+1}=E .
$$

Now let $U$ be a non-zero linear subspace of $\Q^{n+1}$
of dimension $m\leq t_n(\xi )$.
By part (i) of Lemma 7.4,
there is a linear form $L_i\in\{ L_1,L_2\}$ which does not vanish identically
on $U$.
Pick $m-1$ linear forms from $X_0,\ldots ,X_n$
which together with $L_i$ are linearly independent on $U$. Then
the sum of the $H$-exponents corresponding to these linear forms is
$w-m+1$, and again
$$
{w-m+1\over m}\geq {w-{1\over 2}(n+1)+1\over {1\over 2}(n+1)}=E
$$
where we have used $m\leq t_n(\xi )\leq {n+1\over 2}$.
Hence, indeed, the conditions of part (ii) of Lemma 7.3 are satisfied.
This proves part (ii) of Lemma 7.5.
\qed

We now deal with the case that $t_n(\xi )={n+2\over 2}$.
Choose $\mu_0\in\C^*$ such that
$\dim V_n(\mu_0 ,\xi )\hfill\break =t_n(\xi )$
and define
$$
U_0:= \{ {\bf x}\in\Q^{n+1}:\, \sum_{i=0}^n x_iX^i\in V_n(\mu_0 ,\xi )\}
=\{ {\bf x}\in\Q^{n+1}:\, \mu_0\sum_{i=0}^n x_i\xi^i\in \R\}.
\eqno (7.11)
$$
Then $\dim U_0 =t_n(\xi )$ and by Lemma 4.1
the vector space $U_0$ does not depend on the choice of $\mu_0$.
Recall that we can choose $\mu_0$ from $\Q (\xi )$.
Thus, $\mu_0$ is algebraic.

\proclaim Lemma 7.6.
Assume that $t_n(\xi )={n+2\over 2}$ and let $w\geq -1$.
\vskip 2pt\noindent
{\bf (i).} There is a constant $c=c(\xi ,n)>0$ such that for every $H\geq 1$
we have $\lambda_1(\xi ,n,w,H)\leq cH^{{2w-n\over n+2}}$.
\vskip 2pt\noindent
{\bf (ii).}
For every $\vep >0$ there is $H_{2,\vep}>0$  such that for every
$H>H_{2,\vep}$ we have
$$
\eqalignno{
&H^{{2w-n\over n+2}-\vep} < \lambda_1(\xi ,n,w,H)\leq\cdots\leq
\lambda_{(n+2)/2}(\xi ,n,w,H)< H^{{2w-n\over n+2}+\vep}.&(7.12)\cr
&H^{{2w-n+2\over n}-\vep} <
\lambda_{(n+4)/2}(\xi ,n,w,H)\leq\cdots\leq
\lambda_{n+1}(\xi ,n,w,H)< H^{{2w-n+2\over n}+\vep}.
&(7.13)\cr
&H^{{2w-n+2\over n}-\vep}K(\xi ,n,w,H)\cap\Z^{n+1}\subset U_0.
&(7.14)\cr
}
$$

\proof
We first prove part (ii).
The idea is to apply Lemma 7.3 first to a convex body defined on the
quotient space $\R^{n+1}/\R U_0$, and then to $K(\xi ,n,w,H)$ restricted
to $\R U_0$.
   
Let $\mu_0=\alpha_0+\sqrt{-1}\cdot \beta_0$, where
$\alpha_0,\beta_0\in\R$ and define the linear form
$$
M_1:={1\over |\alpha_0|+|\beta_0|}\cdot \Big( \beta_0L_1+\alpha_0L_2\Big).
$$
By a straightforward computation,
$$
M_1={1\over 2 \sqrt{-1}
(|\alpha_0|+|\beta_0|)}\Big( \mu_0\sum_{i=0}^n \xi^iX_i
-\overline{\mu_0}\sum_{i=0}^n \overline{\xi}^iX_i\Big),
$$
hence
$$
\{ {\bf x}\in\Q^{n+1}:\, M_1({\bf x})=0\}=U_0.
\eqno (7.15)
$$
Since $U_0$ has dimension ${n+2\over 2}$,
we can choose linear forms $M_2,\ldots ,M_{n/2}$ in $X_0,\ldots ,X_n$
as follows:
$M_2,\ldots ,M_{n/2}$ vanish identically on $U_0$;
$\{ M_1,M_2,\ldots ,M_{n/2}\}$ is linearly independent;
and each $M_i$ ($i=2,\ldots ,{n\over 2}$)
has real algebraic coefficients the sum of
whose absolute values is equal to $1$.

There is a surjective linear map $\psi$
from $\R^{n+1}$ to $\R^{n/2}$ with kernel $\R U_0$,
which induces a surjective $\Z$-linear map from $\Z^{n+1}$ to $\Z^{n/2}$
with kernel $U_0\cap\Z^{n+1}$.
For $i=1,\ldots ,{n\over 2}$, let $M_i^*$ be the linear form on $\R^{n/2}$
such that $M_i=M_i^*\circ\psi$. Then $M_1^*,\ldots ,M_{n/2}^*$
are linearly independent. Now it is clear that for
${\bf x}\in K(\xi ,n,w,H)$ we have
$$
\eqalign{
&|M_1^*(\psi ({\bf x}))|=|M_1({\bf x})|\leq\max (|L_1({\bf x})|,|L_2({\bf x})|)\leq H^{-w},
\cr
&|M_i^*(\psi ({\bf x}))|=|M_i({\bf x})|\leq\max (|x_0|,\ldots ,|x_n|)\leq H
\ \ (i=2,\ldots , n/2),
}
$$
in other words, if ${\bf x}\in K(\xi ,n,w,H)$ then 
$\psi ({\bf x})$ belongs to the convex body
$$
\Pi (H):=\{ {\bf y}\in \R^{n/2}:\, |M_1^*({\bf y})|\leq H^{-w},\
|M_i^*({\bf y})|\leq H\ (i=2,\ldots ,n/2)\}.
$$
Similarly, for any $\lambda >0$ we have
$$
{\bf x}\in \lambda K(\xi ,n,w,H)\cap\Z^{n+1} \Longrightarrow
\psi ({\bf x})\in \lambda\Pi (H)\cap\Z^{n/2} .
\eqno (7.16)
$$

Let $\vep >0$.
Denote by $\nu_1(H),\ldots ,\nu_{n/2}(H)$ the successive minima of
$\Pi (H)$.
We apply Lemma 7.3.
Let $U$ be a linear subspace of $\Q^{n/2}$
of dimension $m>0$. By (7.15),
$M_1^*$ does not vanish identically on $U$. Pick $m-1$ linear forms
from $M_2^*,\ldots M_{n/2}^*$ which together with $M_1^*$ form a system of
linear forms linearly independent on $U$. The sum of the $H$-exponents
corresponding to these linear forms
is $w-m+1$ and we have
$$
{w-m+1\over m}\geq {2w-n+2\over n}.
$$
So the conditions of part (ii) of Lemma 7.3 are satisfied.
Consequently,
for every sufficiently large $H$ we have
$$
H^{{2w-n+2\over n}-\vep /2n} <\nu_1(H)\leq\cdots\leq \nu_{n/2}(H)
< H^{{2w-n+2\over n}+\vep /2n}.
$$
Together with (7.16) this implies
$$
H^{{2w-n+2\over n}-\vep /2n}K(\xi ,n,w,H)\cap\Z^{n+1}\subset U_0
$$
which implies (7.14).

Further, since $\dim U_0 ={n\over 2}+1$,
we have
$$
H^{{2w-n+2\over n}-(\vep /2n)} <
\lambda_{{n+4\over 2}}(\xi ,n,w,H)\leq\cdots\leq
\lambda_{n+1}(\xi ,n,w,H).
\eqno (7.17)
$$

For $i=1,\ldots ,{n+2\over 2}$,
denote by $\mu_i(H)$ the minimum of all positive
reals $\mu$ such that $\mu K(\xi ,n,w,H)\cap U_0\cap\Z^{n+1}$ contains
$i$ linearly independent points.

We apply again Lemma 7.3.
Let $U$ be a linear subspace of $U_0$ of dimension $m>0$.
By part (i) of Lemma 7.4, there is a linear form $L_i\in\{ L_1,L_2\}$
which does not vanish identically on $U$.
Pick
$m-1$ coordinates from $x_0,\ldots ,x_n$ which together with
$L_i$ form a system of linear forms
which is linearly independent on $U$. Then the sum of the $H$-exponents
corresponding to these linear forms is $w-m+1$ and
$$
{w-m+1\over m}\geq {2w-n\over n+2}.
$$
By means of a bijective linear map $\phi$ from $\R U_0$ to $\R^{(n+2)/2}$
with $\phi (U_0\cap\Z^{n+1} )=\Z^{(n+2)/2}$, we can transform
$K(\xi ,n,w,H)\cap \R U_0$ into a convex body with successive minima
$\mu_1(H),\ldots ,\mu_{(n+2)/2}(H)$ satisfying the conditions of part (ii)
of Lemma 7.3.
It follows that for every sufficiently large $H$,
$$
H^{{2w-n\over n+2}-\vep /2n}< \mu_1(H)\leq\cdots\leq \mu_{n+2\over 2}(H)<
H^{{2w-n\over n+2}+\vep /2n}.
\eqno (7.18)
$$

By combining (7.18) with (7.17) and the already proved (7.14)
we obtain (assuming that $\vep$ is sufficiently small), that
$\mu_i(H)=\lambda_i(\xi ,n,w,H)$ for $i=1,\ldots ,{n+2\over 2}$.
By inserting this into (7.18) we obtain (7.12).

By Minkowski's Theorem,
$$
\prod_{i=1}^{n+1} \lambda_i(\xi ,n,w,H) \ll {\rm Vol}(K(\xi ,n,w,H))^{-1}
\ll H^{2w-n+1\over n+1}.
$$
Together with (7.12), (7.17) this implies that for every
sufficiently large $H$ we have
$$
H^{{2w-n+2\over n}-\vep /2n} <
\lambda_{{n+4\over 2}}(\xi ,n,w,H)\leq\cdots\leq
\lambda_{n+1}(\xi ,n,w,H)< H^{{2w-n+2\over n}+\vep}.
$$
This implies (7.13), and completes the proof of part (ii).

It remains to prove part (i). Applying part (i) of Lemma 7.3
to the image under $\phi$ of $K(\xi ,n,w,H)\cap \R U_0$
we obtain that there is a constant $c=c(\xi ,n)>0$ such that
for every $H\geq 1$ we have $\mu_1(H)\leq H^{{2w-n\over n+2}}$.
Since obviously, $\lambda_1(\xi ,n,w,H)\leq \mu_1(H)$,
part (i) follows.
\qed

\vskip 6mm

\centerline{\bf 8. Proof of Theorem 4}

\vskip 5mm

Let again $n$ be a positive integer, and $\xi$
a complex, non-real algebraic number of degree $>n$.
Let $L_1,L_2$ denote the linear forms defined by (7.8) and
$K(\xi ,n,w,H)$ the convex body defined by (7.9).
Put
$$
u_n(\xi ):= \max \Big\{{n-1\over 2}, t_n(\xi )-1\Big\}.
\eqno (8.1)
$$
In view of Lemma 4.4, in order to prove Theorem 4,
it suffices to prove that $w_n(\xi )\leq u_n(\xi )$,
$\hw_n(\xi )\geq u_n(\xi )$, $w_n^*(\xi )\geq u_n(\xi )$.

\proclaim Lemma 8.1. We have
$w_n(\xi )\leq u_n(\xi )$.

\proof
Let $w\in\R$.
Suppose there are infinitely many polynomials
$P=x_0+x_1X+\cdots +x_nX^n\in\Z[\variable ]$ satisfying
$$
0<|P(\xi )|\leq H(P)^{-w}.
\eqno (8.2)
$$
For such a polynomial $P$, put $H:= H(P)$,
${\bf x}=(x_0,\ldots ,x_n)$.
Then clearly, $|L_1({\bf x})|=|\Re P(\xi )|\leq H^{-w}$,
$|L_2({\bf x})|=|\Im P(\xi )|\leq H^{-w}$,
$|x_i|\leq H$ for $i=0,\ldots ,n$, and so
$$
{\bf x}\in K(\xi ,n,w,H)\cap\Z^{n+1}.
\eqno (8.3)
$$
Since (8.2) is supposed to hold for infinitely many polynomials
$P\in\Z [\variable ]$ of degree $\leq n$, there are arbitrarily large $H$
such that there is a non-zero ${\bf x}$ with (8.3).
That is, there are arbitrarily large $H$ such that the first minimum
$\lambda_1=\lambda_1(\xi ,n,w,H)$ of $K(\xi ,n,w,H)$ is $\leq 1$.

First suppose that $t_n(\xi )\leq {n+1\over 2}$.
Then $u_n(\xi )={n-1\over 2}$. By Lemma 7.5,
for every $\vep >0$ there is $H_{\vep}>1$ such that
$\lambda_1\geq H^{{2w-n+1\over n+1}-\vep}$
for every $H>H_{\vep}$. Hence $w\leq {n-1\over 2}=u_n(\xi )$.
Now suppose that $t_n(\xi )={n+2\over 2}$; then $u_n(\xi )={n\over 2}$.
By Lemma 7.6,
for every $\vep >0$ there is
$H_{\vep}>1$ such that
$\lambda_1\geq H^{{2w-n\over n+2}-\vep}$
for every $H>H_{\vep}$. Hence $w\leq {n\over 2}=u_n(\xi )$.
This implies Lemma 8.1.
\qed

\proclaim Lemma 8.2.
We have $\hw_n(\xi )\geq u_n(\xi )$,
$w_n^*(\xi )\geq u_n(\xi )$.

\proof
We prove the following stronger assertion:
For every $\vep >0$ there is
$H_{\vep}>1$ such that for every $H>H_{\vep}$ there is a non-zero
irreducible polynomial $P\in\Z [\variable ]$ of degree $n$ with
$$
0< |P(\xi )|\leq H^{-u_n(\xi )+\vep},\ \
|P'(\xi )|\geq H^{1-\vep},\ \ H(P)\leq H\, ,
\eqno (8.4)
$$
where $P'$ denotes the derivative of $P$.

By ignoring
the lower bound for $|P'(\xi )|$ in (8.4)
we obtain that for every $H>H_{\vep}$
there is a non-zero irreducible polynomial $P\in\Z [\variable ]$ of degree $n$
such that $0< |P(\xi )|\leq H^{-u_n(\xi )+\vep}$.
This implies $\hw_n (\xi )\geq u_n(\xi )$.

To prove that $w_n^*(\xi )\geq u_n(\xi )$
we have to show that for every $\vep >0$
there are infinitely many algebraic numbers $\alpha$ of degree
at most $n$ with $|\xi -\alpha |\leq H(\alpha )^{-u_n(\xi )-1+ \vep}$.
We prove the existence of infinitely many such $\alpha$ of degree equal
to $n$.  
Take an irreducible polynomial $P\in\Z [X]$ with (8.4)
and let $\alpha$ be a zero of $P$
closest to $\xi$. Then using
the inequalities $|\xi -\alpha |\ll |P(\xi )/P'(\xi )|$
(see (A.11) on p. 228 of \cite{BuLiv}) and $H(\alpha )\ll H(P)\ll H$,
we obtain
$$
|\xi -\alpha | \ll H^{-u_n(\xi )-1+2\vep} 
\ll H(\alpha )^{-u_n(\xi )-1+2\vep},
\eqno (8.5)
$$
where the constants implied by $\ll$ depend only on $n,\vep$.
Since $\deg\xi >n$, the number $\alpha$ cannot be equal to $\xi$ so
eq. (8.5) cannot hold for fixed $\alpha$ and arbitrarily large $H$.
Hence by letting $H\to\infty$, we obtain infinitely many distinct
algebraic numbers $\alpha$ of degree $n$ with (8.5).

We proceed to prove the assertion stated above.
Constants implied by $\ll$ and $\gg$ will depend on $\xi$, $n$ and $\vep$.
Write the polynomial $P$ as
$P=x_0+x_1X+\cdots +x_nX^n$ and put ${\bf x}:=(x_0,\ldots ,x_n)$.
As before, let $L_1,L_2$ be the linear forms given by
$L_1({\bf x})=\Re P(\xi )$, $L_2({\bf x})=\Im P(\xi )$.
Further, define the
linear forms $L_1'$, $L_2'$ by
$$
L_1'({\bf x})=\Re P'(\xi )=\Re (\sum_{j=1}^n j x_j \xi^{j-1}),\quad
L_2'({\bf x})=\Im P'(\xi )=\Im (\sum_{j=1}^n j x_j \xi^{j-1}).
$$
We have to distinguish between the cases $t_n(\xi )\leq {n+1\over 2}$
and $t_n(\xi )={n+2\over 2}$.

First suppose that $t_n(\xi )\leq {n+1\over 2}$.
Then $u_n(\xi )={n-1\over 2}$.
We prove that
for every $\vep >0$ there is $H_{\vep }>1$
with the property that for every $H>H_{\vep }$ there is ${\bf x}\in\Z^{n+1}$
with
$$
\eqalignno{
&|L_1({\bf x})|\leq H^{-{n-1\over 2}+\vep /3},\
|L_2({\bf x})|\leq H^{-{n-1\over 2}+\vep /3},\ \
|x_0|\leq H,\ldots ,|x_n|\leq H&(8.6)\cr
&\max \big\{|L_1'({\bf x})|,|L_2'({\bf x})|\big\}> H^{1-\vep}&(8.7)\cr
&2\not{|} x_n,\ \ 2|x_i\ \hbox{for $i=0,\ldots ,n-1$,}\ \ 4\not{|} x_0.&
(8.8)
}
$$
Then the polynomial $P=\sum_{i=0}^n x_iX^i$ satisfies (8.4)
and by Eisenstein's criterion it is irreducible.

Let $H\geq 1$, $\vep >0$.
Consider vectors
${\bf x}\in\Z^{n+1}$  satisfying (8.6) but not (8.7), i.e., with
$$
\left.
\matrix{
|L_1({\bf x})|\leq H^{-{n-1\over 2}+\vep /3},\
|L_2({\bf x})|\leq H^{-{n-1\over 2}+\vep /3},\hfill\cr
|L_1'({\bf x})|\leq H^{1-\vep},\
|L_2'({\bf x})|\leq H^{1-\vep},\hfill\cr
|x_0|\leq H,\ldots ,|x_n|\leq H.\hfill
}
\right\}
\eqno (8.9)
$$
By considering the coefficients of $X_0,X_1,X_2$ one infers that the linear
forms $L_1$, $L_2$ and $L_2'$ are linearly independent.
Pick $n-2$ coordinates from $X_0,\ldots ,X_n$ which together with
$L_1$, $L_2$, $L_2'$ form a system of $n+1$ linearly independent linear forms.
The sum of the corresponding $H$-exponents is
$$
(n-1-2\vep /3 )+(\vep -1)-(n-2)=\vep /3 >0.
$$
So by Proposition 7.1, there is a finite collection of proper linear
subspaces $T_1,\ldots ,T_m$ of $\Q^{n+1}$, with the property
that for every $H\geq 1$,
there is $T_i\in\{ T_1,\ldots ,T_m\}$ such that the set of solutions
${\bf x}\in\Z^{n+1}$ of (8.9) is contained in $T_i$.
Consequently, if ${\bf x}$ satisfies (8.6) but does not lie in
$T_1\cup\cdots\cup T_m$ then it also satisfies (8.7).

We apply Lemma 7.5 with $w={n-1\over 2}$.
Let $\eta >0$ be a parameter depending on $n,\vep$ to be chosen later,
and $Y$ a parameter depending on $H$ and $\eta$, also chosen later.
For brevity we write $K(Y)$ for the convex body
$K(\xi ,n,{n-1\over 2} ,Y)$ and denote by $\lambda_{n+1}(Y)$
the largest of the successive minima of this body.
According to a result of Mahler (see Lemma 8, page 135 of Cassels \cite{Cas})
there is a constant $c_1=c_1(n)$ such that the convex body
$c_1\lambda_{n+1}(Y)K(Y)$ contains a basis of $\Z^{n+1}$.
By applying Lemma 5.2 with ${\eta\over 2}$ instead of $\vep$
we obtain that for every sufficiently large $Y$ we have
$\lambda_{n+1}(Y)< Y^{\eta/2}$.
Then for every $Y$
large enough to satisfy also $c_1Y^{\eta /2}<Y^{\eta}$,
the
convex body $Y^{\eta}K(Y)$, that is, the body given by
$$
|L_1({\bf x})|\leq Y^{-{n-1\over 2}+\eta},\
|L_2({\bf x})|\leq Y^{-{n-1\over 2}+\eta},\ \
|x_0|\leq Y^{1+\eta},\ldots ,|x_n|\leq Y^{1+\eta}
$$
contains a basis of $\Z^{n+1}$,
$\{ {\bf x}^{(0)},\ldots ,{\bf x}^{(n)}\}$, say.
Consider the vectors
$$
{\bf x}=(x_0,\ldots ,x_n) =\sum_{i=0}^n a_i{\bf x}^{(i)}\
\hbox{with $a_i\in\Z$, $|a_i|\leq Y^{\eta}$ for $i=0,\ldots ,n$.}
\eqno (8.10)
$$
Assuming again that $Y$ is sufficiently large,
each vector (8.10) satisfies
$$
|L_1({\bf x})|\leq Y^{-{n-1\over 2}+3\eta},\
|L_2({\bf x})|\leq Y^{-{n-1\over 2}+3\eta},\ \
|x_0|\leq Y^{1+3\eta},\ldots ,|x_n|\leq Y^{1+3\eta}.
\eqno (8.11)
$$
Since ${\bf x}^{(0)},\ldots ,{\bf x}^{(n)}$ span $\Z^{n+1}$,
the number of vectors (8.10) with the additional property (8.8)
is $\gg Y^{(n+1)\eta}$.
On the other hand,
the number of vectors (8.10) lying in $T_1\cup\cdots\cup T_m$
is $\ll Y^{n\eta}$. Hence if $Y$ is sufficiently large,
there exist vectors ${\bf x}$ with (8.10), (8.8) and with
${\bf x}\not\in T_1\cup\cdots\cup T_m$. Now by
choosing $\eta$ and then $Y$ such that
$$
{{n-1\over 2}-3\eta\over 1+3\eta}={n-1\over 2}-{\vep\over 3},\quad
Y^{1+3\eta}=H\, ,
$$
system
(8.11) translates into (8.6). Thus, we infer that for every sufficiently
large $H$, there exist vectors ${\bf x}\in\Z^{n+1}$ with (8.6), (8.8)
which do not lie in $T_1\cup\cdots\cup T_m$. But as we have seen
above, such vectors satisfy (8.7).
This settles the case that $t_n(\xi )\leq {n+1\over 2}$.
\vskip 5pt

Now assume that $t_n(\xi )={n+2\over 2}$.
Then $u_n(\xi )={n\over 2}$.
We first show that it suffices
to prove that for every $\vep >0$ and every sufficiently large $H$
there exists a polynomial $P\in\Z [\variable ]$ of degree $\leq n$ with (8.4),
without the
requirements that $P$ be irreducible and have degree equal to $n$.
Indeed suppose that for every sufficiently large $H$
there is a polynomial $P\in\Z [\variable ]$
satisfying (8.4) such that $\deg P <n$ or $P$ is reducible.
By the same argument as above, it follows that
there are infinitely many algebraic numbers $\alpha$ of degree $<n$
with (8.5). Then there is $m<n$ such that (8.5) has infinitely
many solutions in algebraic numbers $\alpha$ of degree $m$.
By Lemma 4.2 and our assumption $t_n(\xi )={n+2\over 2}$,
the number $n$ is even, so $n-1$ is odd and hence
$t_{n-1}(\xi )\leq {n\over 2}$. So $u_{n-1}(\xi )={n-2\over 2}<u_n(\xi )$.
Now by Lemmata 4.4 and 8.1,
$$
w_m^*(\xi )\leq w_m(\xi )\leq w_{n-1}(\xi )< u_n(\xi ),
$$
which contradicts that (8.5) has infinitely many solutions in algebraic
numbers $\alpha$ of degree $m$. So for every sufficiently large $H$,
the polynomials $P\in\Z [\variable ]$ of degree $\leq n$ that satisfy (8.4)
necessarily have degree equal to $n$ and are irreducible.

Let $\vep >0$. Let $U_0$ be the vector space defined by (7.11).
Recall that $U_0$ has dimension ${n+2\over 2}$.
We show that for every sufficiently large $H$ there
is a non-zero ${\bf x}\in U_0\cap\Z^{n+1}$ with
$$
\eqalignno{
&|L_1({\bf x})|\leq H^{-{n\over 2}+\vep /3},\  
|L_2({\bf x})|\leq H^{-{n\over 2}+\vep /3},\ \  
|x_0|\leq H,\ldots ,|x_n|\leq H,&(8.12)\cr
&\max \{|L_1'({\bf x})|, |L_2'({\bf x})| \} > H^{1-\vep}.&(8.13)\cr
}
$$

Let $H>1$ and consider those vectors ${\bf x}\in U_0\cap\Z^{n+1}$
satisfying (8.12) but not (8.13), i.e.,
$$
\left.
\matrix{
|L_1({\bf x})|\leq H^{-{n\over 2}+\vep /3},\
|L_2({\bf x})|\leq H^{-{n\over 2}+\vep /3},\hfill\cr
|L_1'({\bf x})|\leq H^{1-\vep},\
|L_2'({\bf x})|\leq H^{1-\vep},\hfill\cr
|x_0|\leq H,\ldots ,|x_n|\leq H.\hfill
}
\right\}
\eqno (8.14)
$$

\noindent
{\bf Claim.}
{\it There are $L_i\in\{ L_1,L_2\}$, $L_j'\in\{ L_1',L_2'\}$
that are linearly independent on $U_0$.}
\vskip 2mm
Assume the contrary.
By Lemma 7.4, the linear forms $L_1,L_2$ are linearly dependent on $U_0$
and at least one of $L_1,L_2$ does not vanish identically on $U_0$.
Hence $M:=L_1+\sqrt{-1}\cdot L_2$ does not vanish identically on $U_0$,
and so $M$ and $M':= L_1'+\sqrt{-1}\cdot L_2'$
are linearly dependent on $U_0$.

Since $\dim U_0 ={n+2\over 2}$, there are two linearly independent vectors
${\bf a}=(a_0,\ldots ,a_n)$,
${\bf b}=(b_0,\ldots ,b_n)\in U_0\cap \Z^{n+1}$
such that if $k$ is the largest
index with $a_k\not= 0$ and $l$ the largest index with $b_l\not= 0$,
then $k<l\leq n-{n+2\over 2}+2={n+2\over 2}$.
Let $A=\sum_{i=0}^k a_iX^i$, $B=\sum_{j=0}^l b_jX^j$.

Then by the linear dependence of $M$, $M'$ we have
$$
A(\xi )B'(\xi )-A'(\xi )B(\xi )=M({\bf a})M'({\bf b})-M'({\bf a})M({\bf b})=0.
$$
But the polynomial $AB'-A'B$ is not identically $0$ 
(since $A,B$ are linearly independent) and has degree at most
$$
k+l-1\leq 2(n+1-(n+2)/2)=n.
$$
This leads to a contradiction since by assumption, $\deg\xi >n$.
This proves our claim.

Choose ${n-2\over 2}$ coordinates from $X_0,\ldots ,X_n$ which together with
$L_i,L_j'$ form a system of ${n+2\over 2}$ linear forms which are linearly
independent on $U_0$. Then the sum of the corresponding $H$-exponents is
$$
({n\over 2}-\vep /3)+(-1+\vep )-{n-2\over 2}= 2\vep /3 >0.
$$
So by Proposition 7.1, there are proper linear subspaces $T_1,\ldots ,T_m$
of $U_0$ with the property that for every $H>1$
there is $T_i\in\{ T_1,\ldots ,T_m\}$ such that the set of
${\bf x}\in U_0\cap\Z^{n+1}$ with (8.14) is contained in $T_i$.
This implies that vectors ${\bf x}\in U_0\cap\Z^{n+1}$ that
satisfy (8.12) and for which ${\bf x}\not\in T_1\cap\cdots\cap T_m$
necessarily have to satisfy (8.13).

To show that there are vectors ${\bf x}\in U_0\cap\Z^{n+1}$ with
${\bf x}\not\in T_1\cup\cdots\cup T_m$ 
one proceeds similarly as above,
but applying Lemma 7.6 with $w={n\over 2}$ instead of Lemma 7.5:
For appropriate $\eta ,Y$, depending on $\vep$, $H$,
one may choose a basis
${\bf x}^{(1)},\ldots ,{\bf x}^{({n+2\over 2})}$
of $U_0\cap\Z^{n+1}$
in $c_2\lambda_{(n+2)/2}(Y)K(Y)$, where $c_2=c_2(n)$ depends only on $n$,
and $\lambda_{(n+2)/2}(Y)$ is the ${n+2\over 2}$-th minimum
of $K(Y):=K(\xi ,n,{n\over 2},Y)$. Then one takes linear combinations
as in (8.10), and by a counting argument
arrives at a vector
${\bf x}$ with (8.12)
which does not lie in $T_1\cup\cdots\cup T_m$, hence satisfies (8.13).
Here, we don't have to impose (8.8).
This completes the proof of Lemma 8.2.
\qed

\vskip 6mm

\centerline{\bf 9. Proof of Theorem 5}

\vskip 5mm

We first prove the following analogue of Theorem 4.

\proclaim Proposition 9.1.
Let $n$ be a positive integer and
$\xi$ a complex, non-real algebraic number of degree $>n$. Then
$$
\tw_n(\xi )=\tw^*_n(\xi )={\hat\tw}_n(\xi )=
n-1-\max\left\{ {n-1\over 2}, t_n(\xi )-1\right\}.
$$

Put $v_n(\xi ):=n-1-\max\{ {n-1\over 2}, t_n(\xi )-1\}$.
Completely similarly as in Lemma 4.4 we have
$$
\tw_n^*(\xi )\leq \tw_n(\xi ),\quad {\hat \tw}_n(\xi )\leq \tw_n(\xi ).
$$
Therefore, in order to prove Proposition 9.1,
it suffices to prove the inequalities
$$
\tw_n^*(\xi) \geq  v_n(\xi ),\ \
{\hat \tw}_n(\xi) \geq v_n(\xi ),\ \
\tw_n(\xi )\leq v_n(\xi ).
$$

These inequalities are proved in Lemmata 9.2 and 9.3 below.
The integer $n$ and the algebraic number $\xi$ will be as in the
statement of Proposition 9.1.

\proclaim Lemma 9.2.
We have
$$
\tw_n^*(\xi) \geq  v_n(\xi ),\quad {\hat \tw}_n(\xi) \geq v_n(\xi ).
$$

\proof
We proceed as in Bugeaud and Teuli\'e \cite{BuTe}, using a method developed
by Davenport and Schmidt \cite{DaSc69}
(see also Theorem 2.11 from \cite{BuLiv}).
As in Section 7, we consider the symmetric convex body
$$
\eqalign{
K(\xi ,n,w,H):= \{ {\bf x}\in\R^{n+1}:\,
|x_n \Re (\xi^n) + \ldots + x_1 \Re (\xi) & + x_0|\leq H^{-w}, \cr
|x_n \Im (\xi^n) + \ldots + x_1 \Im (\xi) & + x_0| \leq H^{-w}, \cr
& |x_0|\leq H,\ldots ,|x_n|\leq H\},\cr}
$$
where ${\bf x}=(x_0,\ldots ,x_n)$ and $w\in\R$.

Set $w:= v_n(\xi )$.
For brevity, we denote the convex body $K(\xi, n ,w,H)$ by $K(H)$.

Let $\vep >0$ be a real number.
Then in both the cases
$t_n(\xi )\leq (n+1)/2$ and $t_n(\xi )=(n+2)/2$ we have,
by Lemmata 7.5 and 7.6, respectively, that
for every sufficiently large real number $H$,
$$
\lambda_{n+1}(H) <H^{\vep},
$$
where $\lambda_{n+1}(H)$ denotes the largest successive mimimum of $K(H)$.

There is a constant $c_1=c_1(n)$ such that the convex body
$c_1\lambda_{n+1}(H)K(H)$ contains a basis of $\Z^{n+1}$,
$$
{\bf x}^{(i)}=(x_0^{(i)},\ldots ,x_n^{(i)})\ \ (i=1,\ldots ,n+1),
$$
say.
This means that there exist $n+1$ integer
polynomials
$$
P_i=x_n^{(i)}X^n+ \ldots +x_1^{(i)}X+x_0^{(i)}, \quad (i=1, \ldots ,n+1),
$$
that form a basis of the $\Z$-module of polynomials in $\Z [\variable ]$ of degree
at most $n$ and for which
$$
H(P_i)\leq c_1H^{1+\vep}, \quad (1\leq i\leq n+1),\eqno (9.1)
$$
and
$$
\max \{|\Re (P_i(\xi))|, |\Im (P_i(\xi) )|\}\leq c_1H^{-w+\vep},
\quad (1\leq i\leq n+1).\eqno (9.2)
$$

There is a unique polynomial $Q=X^{n+1}+\sum_{i=0}^n y_iX^i\in \R [\variable ]$
such that
$$
\left.
\matrix{
\Re Q(\xi )=H^{-w+2\vep},\quad \Im Q(\xi )=H^{-w+2\vep},\quad
\Im Q'(\xi )=H^{1+2\vep},\hfill\cr
y_3=\cdots =y_n=0.
}
\right\}
\eqno (9.3)
$$
Indeed, if we express $\Re Q(\xi )$, $\Im Q(\xi )$ and $\Im Q'(\xi )$
as linear forms in $y_0,\ldots ,y_n$ they form together with
$y_3,\ldots ,y_n$ a linearly independent system of rank $n+1$,
and so (9.3) gives rise to a system of linear equations with a unique
solution $y_0,\ldots ,y_n$.

By expressing
$y_0,y_1,y_2$ as a linear combination of these linear forms, we obtain
$$
|y_i|\ll H^{1+2\vep}\ \ \hbox{for $i=0,1,2$,}
\eqno (9.4)
$$
where here and below, constants implied by $\ll$ depend on $n,\xi ,\vep$
only. Since $P_1,\ldots ,P_{n+1}$ span the vector space of polynomials
with real coefficients of degree at most $n$, there are reals
$\theta_1 ,\ldots ,\theta_{n+1}$ such that
$$
Q =X^{n+1}+2\sum_{i=1}^{n+1}\theta_i P_i .
$$
Now choose integers $t_1,\ldots ,t_n$ with
$$
|\theta_i-t_i|\leq 1\ \ \ (i=1, \ldots , n+1),
\eqno (9.5)
$$
and define the polynomial
$$
P:= X^{n+1}+2\sum_{i=1}^{n+1} t_iP_i .
$$
Write $P=X^{n+1}+\sum_{i=1}^n x_iX^i$.

For a suitable choice of $t_1, \ldots , t_{n+1}$, the polynomial $P$
is irreducible. Indeed,
since $P_1,\ldots ,P_{n+1}$ span the $\Z$-module of all integer polynomials
of degree at most $n$,
at least one of the constant terms $x_0^{(1)},\ldots ,x_0^{(n+1)}$
of $P_1,\ldots ,P_{n+1}$, respectively, must be odd.
Without loss of generality we assume this to be $x_0^{(1)}$.
For a fixed $n$-tuple
$(t_2,\ldots,t_{n+1})$, there are
two choices for $t_1$, that we denote by
$t_{1,0}$ and $t_{1,1}=t_{1,0}+1$.
Since $x_0^{(1)}$ is odd, we can choose $t_1\in\{ t_{1,0},t_{1,1}\}$
such that
$t_1x_0^{(1)}+\cdots+t_{n+1} x_0^{(n+1)}$ is odd.
Then the
constant coefficient of $P(H)$, namely
$2(t_1x_0^{(1)}+\ldots+t_{n+1} x_0^{(n+1)})$, is not divisible by $4$,
and the irreducibility of $P$ follows from the Eisenstein criterion
applied with the prime number $2$.

From (9.5), (9.1), it follows that the absolute values of the
coefficients of $P-Q$ are $\ll H^{1+\vep}$. Further,
by (9.2), (9.1) we have
$$
\eqalign{
&|\Re P(\xi )-\Re Q(\xi )|\ll H^{-w+\vep},\quad
|\Im P(\xi )-\Im Q(\xi )|\ll H^{-w+\vep},\cr
&|\Im P'(\xi )- \Im Q'(\xi )|\ll H^{1+\vep}.
}
$$
Together with (9.3), (9.4) this implies, assuming that
$H$ is sufficiently large,
$$
H(P) \le H^{1 + 3\eps},  \eqno (9.6)
$$
and moreover,
$$
|P(\xi )|\leq |\Re P(\xi )|+|\Im P(\xi )|\leq H^{-w+3\vep},\quad
|P'(\xi )|\geq |\Im P'(\xi )|\geq H^{1+\vep}.
$$

Ignoring the lower bound for $|P'(\xi )|$,
we infer that
$$
{\hat \tw}_n (\xi) \geq (w - 3\eps)/(1 + 3\eps).
$$
Since $\eps$ is arbitrary, we get the second statement of the lemma.
Furthermore, we deduce that the monic polynomial
$P$ has a complex root $\alpha$ with
$$
|\xi - \alpha| \ll {|P(\xi)| \over |P'(\xi)|}
\ll \, H(\alpha)^{- (w + 1 - 2\eps)/(1 + 3 \eps)}.
$$
Since $\eps$ is arbitrary, this shows that
$$
\tw_n^*(\xi) \geq w = v_n(\xi ),
$$
and the proof of Lemma 9.2 is complete. \cqfd

We now prove an upper bound for $\tw_n(\xi )$.

\proclaim Lemma 9.3.
We have
$$
\tw_n(\xi )\leq v_n(\xi ).
\eqno (9.7)
$$

\proof
It suffices to show that for every $w>v_n(\xi )$, the inequality
$$
0<|P(\xi )|\leq H(P)^{-w}
\eqno (9.8)
$$
has only finitely many solutions in monic polynomials $P\in \Z [\variable ]$
of degree at most $n+1$. By replacing any monic polynomial $P$
of degree $k<n+1$ satisfying (9.8) by $X^{n-k}P$
and modifying $w$ a little bit,
one easily observes that it suffices
to show that for every $w>v_n(\xi )$, inequality (9.8)
has only finitely many solutions in monic polynomials $P\in \Z [\variable ]$
of degree {\it precisely} $n+1$.

We have again to distinguish between the cases $t_n(\xi )\leq {n+1\over 2}$
and $t_n(\xi )={n+2\over 2}$. The first case is dealt with
by a modification of the proof of Lemma 7.5,
and the second by a modification of the proof of Lemma 7.6.

First assume that $t_n(\xi )\leq {n+1\over 2}$.
Then $v_n(\xi )={n-1\over 2}$.
Consider the inequality (9.8)
to be solved in monic polynomials $P\in \Z [\variable ]$ of degree $n+1$.
Defne the polynomial $P=\sum_{i=0}^{n+1} x_iX^i$ where $x_{n+1}=1$ and put
${\bf x}=(x_0,\ldots ,x_n,x_{n+1})$, $H:= H(P)$. Define the linear
forms
$$
\tilde{L}_1:=\Re\big( \sum_{i=0}^{n+1}\xi^iX_i\big),\ \
\tilde{L}_2:=\Im\big( \sum_{i=0}^{n+1}\xi^iX_i\big),\ \
\tilde{M}:= \sum_{i=0}^{n+1}\xi^iX_i\, .
$$
Then we can translate (9.8) into
$$
\left.
\matrix{
|\tilde{L}_1({\bf x})|\leq H^{-w},\ |\tilde{L}_2({\bf x})|\leq H^{-w},\hfill\cr
|x_0|\leq H ,\ldots ,|x_n|\leq H,\ |x_{n+1}|\leq 1,\ x_{n+1}\not= 0,\
\tilde{M}({\bf x})\not= 0.
}
\right\}
\eqno (9.9)
$$
We prove that for every $w>{n-1\over 2}$ there is $H_w>1$ such that
if (9.9) has a solution ${\bf x}\in\Z^{n+2}$ then $H<H_w$.
This implies at once that for every $w>{n-1\over 2}$
there are only finitely many monic polynomials $P\in \Z [\variable ]$
of degree $\leq n+1$ with (9.8), and hence that
$\tw_n(\xi )\leq {n-1\over 2}=v_n(\xi )$.

We apply Lemma 7.2. Let $w>{n-1\over 2}$. We have to verify (7.3).
First, let $U$ be a linear subspace of $\Q^{n+2}$
of dimension $m>{n+3\over 2}$
on which $X_{n+1}$ and $\tilde{M}$ are not identically $0$.
Then $\tilde{L}_1$, $\tilde{L}_2$, $X_{n+1}$ are linearly independent on $U$.
For if not, then the linear forms
$$
L_1:=\Re\big( \sum_{i=0}^n\xi^iX_i\big),\ \
L_2:=\Im\big( \sum_{i=0}^n \xi^iX_i\big)
$$
are linearly dependent on $U\cap \{x_{n+1}=0\}$ 
which has dimension
larger than ${n+1\over 2}$. But by part (ii) of Lemma 7.4 this is impossible.
Now choose $m-3$ coordinates from $X_0,\ldots, X_n$ which together with
$\tilde{L}_1$, $\tilde{L}_2$, $X_{n+1}$ are linearly independent on $U$.
Then the $H$-exponents corresponding to these linear forms have sum
$$
2w+0+(3-m)> n+2-m\geq 0.
$$

Now let $U$ be a linear subspace of $\Q^{n+2}$ of dimension $m$ with
$2\leq m\leq {n+3\over 2}$ on which $X_{n+1}$, $\tilde{M}$
do not vanish identically.
Then there is $\tilde{L}_i\in\{ \tilde{L}_1,\tilde{L}_2\}$ such that
$\tilde{L}_i$ and $X_{n+1}$ are linearly independent on $U$.
For if not then both $L_1$ and $L_2$ vanish identically on
$U\cap\{ x_{n+1}=0\}$ 
which is impossible by part (i) of Lemma 7.4.
Choose $m-2$ coordinates from $X_0,\ldots ,X_n$ which together with
$\tilde{L}_i$ and $X_{n+1}$ are linearly independent on $U$. Then the
$H$-exponents corresponding to these linear forms have sum
$$
w+0+(2-m) > {n-1\over 2}+2-m\geq 0.
$$

Finally, let $U$ be a one-dimensional linear subspace of $\Q^{n+2}$ on
which none of $X_{n+1}$, $\tilde{M}$, vanishes identically.
Then there is $\tilde{L}_i\in\{ \tilde{L}_1,\tilde{L}_2\}$
not vanishing identically on $U$, and the $H$-exponent corresponding to
this linear form is $w>0$. We conclude that condition (7.3) of Lemma 7.2
is satisfied. So indeed there is $H_w>0$ such that if (9.9) is satisfied
by some ${\bf x}\in\Z^{n+1}$ then $H<H_w$. This settles the case
that $t_n(\xi )\leq {n+1\over 2}$.
\vskip 2pt

Now assume that $t_n(\xi )={n+2\over 2}$.
Then $v_n(\xi )={n-2\over 2}$. Further,
by Lemmata 4.2 and 4.3,
$n$ is even, $n+1<\deg\xi$, and
$$
t_{n+1}(\xi )=t_n(\xi )={n+2\over 2}.
\eqno (9.10)
$$
Choose $\mu_0=\alpha_0+\sqrt{-1}\cdot \beta_0$ with $\alpha_0,\beta_0\in\R$
such that $\dim V_n(\mu_0,\xi )=t_n(\xi )={n+2\over 2}$.
Define the linear form
$$
\tilde{M}_1={1\over |\alpha_0|+|\beta_0|}\cdot
\Big( \beta_0\tilde{L}_1+\alpha_0\tilde{L}_2\Big)
={1\over 2 \sqrt{-1}
(|\alpha_0|+|\beta_0|)}\Big( \mu_0\sum_{i=0}^{n+1} x_i\xi^i
-\overline{\mu_0}\sum_{i=0}^{n+1} x_i\overline{\xi}^i\Big).
$$
Let
$$
\tilde{U}_0=\{ {\bf x}\in\Q^{n+2}:\, \tilde{M}_1({\bf x})=0\}.
$$
Then ${\bf x}=(x_0,\ldots ,x_{n+1})\in \tilde{U}_0$ if and only if
$\sum_{i=0}^{n+1} x_iX^i\in V_{n+1}(\mu_0,\xi )$.

We claim that $X_{n+1}=0$ identically on $\tilde{U}_0$.
Suppose $\tilde{U}_0$ contains a vector ${\bf x}=(x_0,\ldots ,x_{n+1})$
with $x_{n+1}\not= 0$. Then the polynomial
$\sum_{i=0}^{n+1} x_iX^i$ belongs to
$V_{n+1}(\mu_0 ,\xi )$ but not to $V_n(\mu_0 ,\xi )$ which is
impossible by (9.10).
This argument shows also that
$\dim \tilde{U}_0=\dim V_n(\mu_0,\xi )={n+2\over 2}$.

There are linear forms $\tilde{M}_2,\ldots ,\tilde{M}_{n/2}$ in
$X_0,\ldots ,X_{n+1}$ with the following properties:
$\tilde{M}_2,\ldots ,\tilde{M}_{n/2}$ vanish indentically
on $\tilde{U}_0$;
$\{ \tilde{M}_1,\tilde{M}_2,\ldots ,\tilde{M}_{n/2},X_{n+1}\}$
is linearly independent;
and each $\tilde{M}_i$ $(i=2,\ldots ,{n\over 2})$
has real algebraic coefficients
whose absolute values have sum equal to $1$.

Let $\psi$ be a surjective linear mapping from
$\R^{n+2}$ to $\R^{{n+2\over 2}}$
with kernel $\R \tilde{U}_0$ such that the restriction of $\psi$ to $\Z^{n+2}$
maps surjectively to $\Z^{{n+2\over 2}}$
and has kernel $\tilde{U}_0\cap\Z^{n+2}$.
For $i=1,\ldots {n\over 2}$,
let $\tilde{M}_i^*$ be the linear form on $\R^{{n+2\over 2}}$
with $\tilde{M}_i=\tilde{M}_i^*\circ \psi$. Further, let
$\tilde{M}_0^*$ be the linear form on $\R^{{n+2\over 2}}$ such that
$X_{n+1}=\tilde{M}_0^*\circ\psi$.
Then $\tilde{M}_0^*,\ldots ,\tilde{M}_{n/2}^*$ are linearly independent.

Let $w> v_n(\xi )={n-2\over 2}$.
Let $P\in\Z [\variable ]$ be a monic polynomial of degree $n+1$ satisfying (9.8).
Write $P=\sum_{i=0}^{n+1} x_iX^i$, $x_{n+1}=1$,
${\bf x}=(x_0,\ldots ,x_{n+1})$, $H:= H(P)$.
Then ${\bf x}$ satisfies (9.9). By an easy computation it follows
that ${\bf y}:=\psi ({\bf x})$ satisfies
$$
\left.
\matrix{
|\tilde{M}_1^*({\bf y})|\leq H^{-w},\
|\tilde{M}_i^*({\bf y})|\leq H\ (i=2,\ldots ,n/2),\hfill\cr
|\tilde{M}_0^*({\bf y})|\leq 1,\ M_0^*({\bf y})\not= 0.\hfill
}
\right\}
\eqno (9.11)
$$
We show that system (9.11) satisfies condition (7.3) of Lemma 7.2.
First let $U=\Q^{{n+2\over 2}}$.
As observed before, the linear forms
$\tilde{M}_0^*,\ldots ,\tilde{M}_{n/2}^*$ are linearly independent,
and the $H$-exponents corresponding to these linear forms have sum
$$
w-(n-(n/2)-1)+0 >0.
$$
Now let $U$ be a linear subspace of $\Q^{{n+2\over 2}}$
of dimension $m$ with
$0<m\leq {n\over 2}$ on
which $\tilde{M}_0^*$ does not vanish identically.
The linear form $\tilde{M}_1$ can not vanish identically on
$\psi^{-1}(U)$ since $\psi^{-1}(U)$ is strictly larger than $U_0$,
therefore, $\tilde{M}_1^*$ does not vanish identically on $U$.
Choose $m-1$ linear forms among
$\tilde{M}_0^*,\tilde{M}_2^*,\ldots ,\tilde{M}_{n/2}^*$
which together with $\tilde{M}_1^*$ are linearly independent on $U$.
Then the sum of the $H$-exponents corresponding to these linear forms
is at least
$$
w-(m-1)\geq w-((n/2)-1)>0.
$$
Hence condition (7.3) of Lemma 7.2 is satisfied. It follows that there
is $H_w>0$ such that if system (9.11)
is solvable in ${\bf y}\in\Z^{{n+2\over 2}}$
then $H\leq H_w$. Hence for every monic polynomial $P\in \Z [\variable ]$ of degree
$n+1$ with (9.8) we have $H(P)\leq H_w$, implying that (9.8) has only
finitely many solutions.
\qed

As observed above, Proposition 9.1 follows from Lemmata 9.2 and 9.3.\qed

\noindent
{\bf Proof of Theorem 5.} We first prove (3.1).
Assume that $\deg\xi =:d\leq n+1$.
By Liouville's inequality (4.1) we have $\tw_n(\xi )\leq {d-2\over 2}$.
By Proposition 9.1 and Lemma 4.3 we have
$\tw_{d-1}^*(\xi )={\hat\tw}_{d-1}(\xi )={d-2\over 2}$. Hence
$$
\tw_n^*(\xi )\geq \tw_{d-1}^*(\xi )={d-2\over 2},\ \
{\hat\tw}_n(\xi )\geq {\hat\tw}_{d-1}(\xi )={d-2\over 2}.
$$
These facts together imply (3.1).

The equalities (3.2) and (3.3) follow at once by combining
Proposition 9.1 with part (ii) of Lemma 4.2.
The last assertion of Theorem 5 follows at once from Theorem 4
and Proposition 9.1. This completes the proof of Theorem 5.
\qed

\vskip 6mm

\centerline{\bf 10. A refined question}

\vskip 5mm

The exponents $w_n, {\hat w}_n, \ldots $ are defined as suprema of certain
sets of real numbers. We may further ask whether the suprema are
also maxima. In other words, for a given complex number $\xi$,
a positive integer $n$, do there exist a constant $c(\xi, n)$
and infinitely many integer polynomials $P(H)$ of degree at most $n$
such that
$$
0 < |P(\xi)| \le c(\xi, n) \, H(P)^{-w_n(\xi)} \ ?
$$
This is Problem P.1, page 210, of \cite{BuLiv}.

When $\xi$ is algebraic and
real, the answer is clearly positive, by Dirichlet's Theorem.
When $\xi$ is algebraic and non-real, we have already seen that
$w_n(\xi)$ can be much larger than expected; however, the answer to the above
question is also positive.

\proclaim Proposition 10.1.
For any positive integer $n$ and any complex, non-real algebraic number $\xi$,
there exist a constant $c(\xi, n)>0$
and infinitely many integer polynomials $P(H)$ of degree at most $n$
such that
$$
0 < |P(\xi)| \le c(\xi, n)\, H(P)^{-w_n(\xi)}.
\eqno (10.1)
$$

\noi {\bf Proof: }
This follows from (the proof of) Satz 1 from Schmidt \cite{SchmMA}; however,
we feel that it is better to include a complete proof.
Constants implied by $\ll$, $\gg$ depend only on $n,\xi$.

First assume that $d:= \deg \xi >n$. We apply part (i)
of Lemmata 7.5 and 7.6, respectively, with
$w=w_n(\xi )$. Then in view of Theorem 4,
in both the cases $t_n(\xi )\leq {n+1\over 2}$, $t_n(\xi )={n+2\over 2}$,
we have that for every $H\geq 1$
the first minimum $\lambda_1(\xi ,n,w,H)$ of the convex body
$K(\xi ,n,w,H)$ defined by (7.9) is $\ll 1$. Consequently, for
every $H\geq 1$, there is a non-zero polynomial
$P=\sum_{i=0}^n x_i\variable^i\in\Z [\variable ]$
such that
$$
\eqalign{
&|\Re P(\xi )|=|L_1({\bf x})|\ll H^{-w},\
|\Im P(\xi )|=|L_2({\bf x})|\ll H^{-w},\cr
&H(P)=\max \{|x_0|,\ldots ,|x_n|\}\ll H.\cr
}
$$
This clearly implies $|P(\xi )|\ll H^{-w}\ll H(P)^{-w}$.
Arbitrarily large $H$ cannot give rise to the same polynomial $P$
since otherwise we would have $P(\xi )=0$,
against our assumption that $\deg \xi >n$. This proves Proposition 10.1
in the case that $d>n$.

To treat the case $n\geq d$ we simply have to observe
that by Theorem 1
we have $w_n(\xi )=w_{d-1}(\xi )={d-2\over 2}$ and that by what we have proved
above, (10.1) has already infinitely many solutions in polynomials
$P$ of degree at most $d-1$.
\qed

Actually, the above proof yields that the analogue of Proposition 10.1 is
true for the uniform exponent of approximation ${\hat w}_n$.
However, it is a very interesting, but presumably very difficult,
question to decide whether the analogue
of Proposition 10.1 holds for the exponent $w_n^*$.

We briefly summarize what is known on this question.

\proclaim Proposition 10.2.
For any positive integer $n$ and any complex algebraic number $\xi$
of degree $n+1$,
there exist a constant $c(\xi)$
and infinitely many algebraic numbers $\alpha$ of degree at most $n$
such that
$$
0 < |\xi - \alpha| \le c(\xi) \, H(\alpha)^{-w^*_n(\xi)-1}.
$$

\noi {\bf Proof: }
When $\xi$ is real, Proposition 10.2 has been established by
Wirsing \cite{Wir} (see also Theorem 2.9 in \cite{BuLiv}, which reproduces
an alternative proof, due to Bombieri and Mueller \cite{BoMu}).
Without any additional complication, the same method gives the required
result when $\xi$ is complex and non-real.\cqfd

Furthermore, Davenport and Schmidt \cite{DaSc67} proved that
for every real algebraic number $\xi$ of degree at least $3$, 
there exist a constant $c(\xi)$
and infinitely many algebraic integer $\alpha$ of degree at most $2$
such that
$$
0 < |\xi - \alpha| \le c(\xi) \, H(\alpha)^{-w^*_2(\xi)-1}
= c(\xi) \, H(\alpha)^{-3}.
$$
This is a consequence of a more general result
of theirs on linear forms \cite{DaSc68,DaSc70},
which is the key tool for the proof of
the second assertion of the next proposition.

\proclaim Proposition 10.3.
{\bf (i).} For any complex algebraic number $\xi$ of degree greater than $2$,
there exist a constant $c(\xi)$
and infinitely many algebraic numbers $\alpha$ of degree at most $2$
such that
$$
0 < |\xi - \alpha| \le c(\xi) \, H(\alpha)^{-w^*_2(\xi)-1}.
$$
{\bf (ii).} For any complex algebraic number $\xi$ of degree greater than $4$
satisfying $w^*_4(\xi) = 2$,
there exist a constant $c(\xi)$
and infinitely many algebraic numbers $\alpha$ of degree at most $4$
such that
$$
0 < |\xi - \alpha| \le c(\xi) \, H(\alpha)^{-w^*_4(\xi)-1}.
$$

\noi {\bf Proof: }
Let $\xi$ be a complex non-real number of degree greater than $2$.
By the proof of Proposition 10.1,
there are infinitely many integer quadratic polynomials $P$ satisfying
$$
|P(\xi)| \ll H(P)^{-w_2(\xi)}, \quad |P(\ovxi)| \ll H(P)^{-w_2(\xi)}.
$$
Such a polynomial $P$ has a root very near to $\xi$ and another
very near to $\ovxi$. Consequently, it satisfies $|P'(\xi)| \gg H(P)$ and
its root $\alpha$ near to
$\xi$ is such that $|\xi - \alpha| \ll H(\alpha)^{-w_2(\xi)-1}$.
This proves the first part of the proposition since $w_2(\xi) = w_2^*(\xi)$.

Let $\xi$ be a complex (non-real) algebraic number
of degree $>4$ satisfying $w^*_4(\xi) = 2$.
By Theorem 4, this means that $t_4(\xi )=3$,
i.e., there is $\mu_0$ such that $\dim V_4(\mu_0 ,\xi )= 3$.
Let $U_0$ be the vector space of ${\bf x}=(x_0,\ldots ,x_4)\in\Q^5$
such that $\sum_{i=0}^4 x_iX^i\in V_4(\mu_0,\xi )$. Define the linear
forms $L_1,L_2,L_1' ,L_2'$ by
$$
L_1({\bf x})= \Re P(\xi ),\ L_2({\bf x})=\Im P(\xi ),\
L_1'({\bf x})=\Re P'(\xi ),\ L_2'({\bf x})=\Im P'(\xi ),
$$
where $P=\sum_{i=0}^4 x_iX^i$. By Lemma 7.4, the linear forms $L_1,L_2$
are linearly dependent on $U_0$. On the other hand,
by the Claim in the proof of Lemma 8.2,
there are $i,j\in\{ 1,2\}$ such that $L_i,L_j'$ are linearly independent on
$U_0$. Choose linearly independent polynomials $P_1,P_2,P_3$ from
$V_4(\mu_0,\xi )$ with integer coefficients.
By Lemma 4.2 we may assume $\deg P_1<\deg P_2<\deg P_3=4$.
Express $P\in V_4(\mu_0,\xi )$
as $y_1P_1+y_2P_2+y_3P_3$ with ${\bf y}=(y_1,y_2,y_3)\in\Q^3$.
Thus, $\Re P(\xi )$, $\Im P(\xi )$, $\Re P'(\xi )$, $\Im P'(\xi )$
can be expressed as linear forms in ${\bf y}$,
$$
\Re P(\xi )= M_1({\bf y}),\ \Im P(\xi )= M_2({\bf y}),\
\Re P'(\xi )= M_1'({\bf y}),\ \Im P'(\xi )= M_2'({\bf y})
$$
say,
and by the above, $M_1,M_2$ are linearly dependent and there are
$i,j\in\{ 1,2\}$ such that $M_i,M_j'$ are linearly independent.

By Theorem 1 from \cite{DaSc68},
there are infinitely many integer triples ${\bf y}=(y_1,y_2,y_3)$ with
$$
|M_i({\bf y})| \ll |M_j'({\bf y})|\times \| {\bf y}\|^{-3},
$$
where $\|{\bf y}\|=\max\{|y_1|, |y_2|, |y_3|\}$.
This implies that there are infinitely many integer polynomials $P$
of degree $4$
of the shape $y_1P_1+y_2P_2+y_3P_3$ with $y_1,y_2,y_3\in\Z$
such that
$$
{|P(\xi)| \over |P'(\xi)|} \ll H(P)^{-3}.
$$
Consequently, there are infinitely many algebraic numbers $\alpha$
of degree at most $4$ such that $|\xi - \alpha| \ll H(\alpha)^{-3}$.
This completes the proof of Proposition 10.3. \cqfd

\vskip 7mm
\noi {\bf Acknowledgements.}
We are pleased to thank
Noriko Hirata-Kohno, Corentin Pontreau and
Damien Roy for helpful discussions.
The research leading to this paper started with a discussion by both authors
at the Erwin Schroedinger Institut in April 2006,
in the frame of a research program on Diophantine
approximation and heights organized by David Masser, Hans Peter
Schlickewei and Wolfgang Schmidt.
\vskip 8mm

\centerline{\bf References}

\vskip 7mm

\beginthebibliography{999}

\bibitem{BoMu}
E. Bombieri and J. Mueller,
{\it Remarks on the approximation to an
algebraic number by algebraic numbers}, Mich. Math. J.
{33} (1986), 83--93.

\bibitem{BuLiv}
Y. Bugeaud,
Approximation by algebraic numbers.
Cambridge Tracts in Mathematics 160, Cambridge University Press, 2004.

\bibitem{BuLau}
Y. Bugeaud and M. Laurent,
{\it Exponents of Diophantine approximation and
Sturmian continued fractions},
Ann. Inst. Fourier (Grenoble) 55 (2005), 773--804.

\bibitem{BuLauB}
Y. Bugeaud and M. Laurent,
{\it Exponents of homogeneous and inhomogeneous Diophantine
Approximation},
Moscow Math. J. 5 (2005), 747--766.

\bibitem{BuTe}
Y. Bugeaud et O. Teuli\'e,
{\it Approximation d'un nombre r\'eel par des nombres
alg\'e\-briques de degr\'e donn\'e},
Acta Arith. 93 (2000), 77--86.

\bibitem{Cas}
J. W. S. Cassels,
An Introduction to the Geometry of Numbers.
Springer Verlag, 1997.

\bibitem{DaSc67}
H. Davenport and W. M. Schmidt,
{\it Approximation to real numbers
by quadratic irrationals}, Acta Arith. 13 (1967), 169--176.

\bibitem{DaSc68}
H. Davenport and W. M. Schmidt,
{\it A theorem on linear forms},
Acta Arith. 14 (1967/1968), 209--223.

\bibitem{DaSc69}
{H. Davenport and W. M. Schmidt},
{\it Approximation to real numbers by
algebraic integers}, Acta Arith. {15} (1969), 393--416.

\bibitem{DaSc70}
{H. Davenport and W. M. Schmidt},
{\it Dirichlet's theorem on
Diophantine approximation. II}, Acta Arith. {16} (1970), 413--423.

%

\bibitem{EvSchl}
J.-H. Evertse and H.P. Schlickewei,
{\it A quantitative version of the Absolute Parametric Subspace Theorem},
J. reine angew. Math. 548 (2002), 21-127.




\bibitem{Kh}
A. Ya. Khintchine,
{\it \"Uber eine Klasse linearer diophantischer Approximationen},
Rendiconti Circ. Mat. Palermo 50 (1926), 170--195.

\bibitem{Kok}
J. F. Koksma,
{\it \"Uber die Mahlersche Klasseneinteilung der transzendenten Zahlen
und die Approximation komplexer Zahlen durch algebraische Zahlen},
Monatsh. Math. Phys. 48 (1939), 176--189.

\bibitem{Mah}
K. Mahler,
{\it Zur Approximation der Exponentialfunktionen und des
Logarithmus. I, II},
J. reine angew. Math. 166 (1932), 118--150.


\bibitem{Ro}
{K. F. Roth},
{\it Rational approximations to algebraic numbers},
{Mathematika} {2} (1955), 1--20; corrigendum, 168.

\bibitem{RoyA}
D. Roy,
{\it Approximation simultan\'ee d'un nombre et son carr\'e},
C. R. Acad. Sci. Paris 336 (2003), 1--6.

\bibitem{RoyB}
D. Roy,
{\it Approximation to real numbers by cubic algebraic numbers, I},
Proc. London Math. Soc. 88 (2004), 42--62.

\bibitem{RoyC}
D. Roy,
{\it Approximation to real numbers by cubic algebraic numbers, II},
Ann. of Math. 158 (2003), 1081--1087.

\bibitem{RoyWa}
D. Roy and M. Waldschmidt,
{\it Diophantine approximation by
conjugate algebraic integers},
Compos. Math. 140 (2004), 593--612.

\bibitem{SchmAM}
{W. M. Schmidt},
{\it Simultaneous approximation to algebraic numbers by rationals},
Acta Math. 125 (1970), 189--201.

\bibitem{SchmMA}
W. M. Schmidt,
{\it Linearformen mit algebraischen Koeffizienten. II},
Math. Ann. 191 (1971), 1--20.

\bibitem{SchmEM}
W. M. Schmidt,
Approximation to algebraic numbers. Monographie
de l'En\-seigne\-ment Math\'ematique {19}, Gen\`eve, 1971.

\bibitem{SchmLN}
{W. M. Schmidt},
Diophantine Approximation. Lecture Notes in Math.
{785}, Springer, Berlin, 1980.

\bibitem{Spr}
{V. G. Sprind\v zuk},
Mahler's problem in metric number theory.
Izdat. \og Nauka i Tehni\-ka \fg , Minsk, 1967 (in  Russian).
English translation by B. Volkmann, Translations of Mathematical Monographs,
Vol. 25, American Mathematical Society, Providence, R.I., 1969.


\bibitem{Wir}
{E. Wirsing},
{\it Approximation mit algebraischen Zahlen beschr\"ankten
Grades}, J. reine angew. Math. {206} (1961), 67--77.

\endthebibliography

\vskip1cm
\vbox{
\hbox{Yann Bugeaud \hfill}
\hbox{Universit\'e Louis Pasteur \hfill}
\hbox{U. F. R. de math\'ematiques \hfill}
\hbox{7, rue Ren\'e Descartes \hfill}
\hbox{67084 STRASBOURG (FRANCE)\hfill}
\vskip2pt\noindent
\hbox{{\tt bugeaud@math.u-strasbg.fr} \hfill}
}
\vskip1cm\noindent
\vbox{
\hbox{Jan-Hendrik Evertse \hfill}
\hbox{Universiteit Leiden \hfill}
\hbox{Mathematisch Instituut \hfill}
\hbox{Postbus 9512 \hfill}
\hbox{2300 RA LEIDEN (THE NETHERLANDS)\hfill}
\vskip2pt\noindent
\hbox{{\tt evertse@math.leidenuniv.nl} \hfill}
}

\bye